\documentclass[12pt,reqno,a4wide]{article}
\usepackage[notref,notcite]{showkeys}
\usepackage{eurosym}
\usepackage{a4wide}
\usepackage{comment}
\usepackage{mathptmx}
\usepackage{amssymb}
\usepackage{amsfonts}
\usepackage{amsthm}
\usepackage{amsmath}
\usepackage{dsfont}
\usepackage{graphicx}
\usepackage{shadow}
\usepackage[all]{xy}
\usepackage{enumerate}
\usepackage{makeidx}
\usepackage{mathrsfs}
\usepackage{upgreek}
\usepackage{stmaryrd}
\usepackage[utf8]{inputenc}
\usepackage{array}
\usepackage{tikz-cd}
\usepackage{thumbpdf}
\usepackage[colorlinks=true,pagebackref]{hyperref}
\usepackage{textcomp}
\usepackage{chapterbib}
\usepackage{graphics}
\usepackage{longtable}
\usepackage{epsf}
\usepackage{bm}
\usepackage{adjustbox}
\usepackage{centernot}

\makeindex
\newtheorem{theorem}{Theorem}[section]
\newtheorem{lem}[theorem]{Lemma}
\newtheorem{thm}[theorem]{Theorem}
\newtheorem{prop}[theorem]{Proposition}
\newtheorem{cor}[theorem]{Corollary}
\theoremstyle{definition}
\newtheorem*{Beweis}{Proof}

\newtheorem{definition}[theorem]{Definition}
\newtheorem{defns}[theorem]{Definitions}
\newtheorem{rem}[theorem]{Remark}
\newtheorem{rems}[theorem]{Remarks}
\newtheorem{punto}[theorem]{}

\theoremstyle{remark}

\newtheorem{ex}[theorem]{Example}
\newtheorem{exs}[theorem]{Examples}

\input{tcilatex}
\begin{document}

\title{Pushouts and e-Projective Semimodules\thanks{%
MSC2010: Primary 18G05; Secondary 18A30, 16Y60 \newline
Key Words: Semirings; Semimodules; Pushouts Projective Semimodules; Exact
Sequences \newline
The authors would like to acknowledge the support provided by the Deanship
of Scientific Research (DSR) at King Fahd University of Petroleum $\&$
Minerals (KFUPM) for funding this work through projects No. RG1304-1 $\&$
RG1304-2}}
\author{$%
\begin{array}{ccc}
\text{Jawad Abuhlail}\thanks{\text{Corresponding Author}} &  & \text{Rangga
Ganzar Noegraha}\thanks{\text{The paper is extracted from his Ph.D.
dissertation under the supervision of Prof. Jawad Abuhlail.}} \\ 
\text{abuhlail@kfupm.edu.sa} &  & \text{rangga.gn@universitaspertamina.ac.id}
\\ 
\text{Department of Mathematics and Statistics} &  & \text{Universitas
Pertamina} \\ 
\text{King Fahd University of Petroleum $\&$ Minerals} &  & \text{Jl. Teuku
Nyak Arief} \\ 
\text{31261 Dhahran, KSA} &  & \text{Jakarta 12220, Indonesia}%
\end{array}%
$}
\date{\today }
\maketitle

\begin{abstract}
Projective modules play an important role in the study of the category of
modules over rings and in the characterization of various classes of rings.
Several characterizations of projective objects which are equivalent for
modules over rings are \emph{not} necessarily equivalent for semimodules
over an arbitrary semiring. We study several of these notions, in particular
the $e$\emph{-projective semimodules } introduced by the first author using
his new notion of \textit{exact sequences} of semimodules. As pushouts of
semimodules play an important role in some of our proofs, we investigate
them and give a constructive proof of their existence in a way that proved
be very helpful.
\end{abstract}


\section*{Introduction}

The importance of \emph{semirings} (defined, roughly, as rings not
necessarily with subtraction) stems from the fact that they can be
considered as a generalization of both rings and distributive bounded
lattices. Moreover, {semirings, and their \emph{semimodules} (defined,
roughly, as modules not necessarily with subtraction), proved to have wide
applications in many aspects of Computer Science and Mathematics, e.g.,
Automata Theory \cite{HW1998}, Tropical Geometry \cite{Gla2002} and
Idempotent Analysis \cite{LM2005}. Many of these applications can be found
in Golan's book \cite{Gol1999}, which is our main reference in this topic. }

\bigskip

The notion of projective objects can be defined in any category relative to
a suitable \emph{factorization system} of its arrows. Projective semimodules
have been studied intensively (see \cite{Gla2002} for details). {Recently, }%
several papers by Abuhlail, I'llin, Katsov and Nam (among others) prepared
the stage for a homological characterization of special classes of semirings
using special classes of projective, injective and flat semimodules (cf., 
\cite{KNT2009}, \cite{Ili2010}, \cite{KN2011}, \cite{Abu2014}, \cite{KNZ2014}%
, \cite{AIKN2015}, \cite{IKN2017}, \cite{AIKN2018}). For example, \emph{%
ideal-semisimple} semirings all of whose left cyclic semimodules are
projective have been investigated in \cite{IKN2017}.

\bigskip

In addition to the \emph{categorical notions} of \emph{projective semimodules%
} over a semiring, several other notions were considered in the literature,
e.g., the so called $k$-projective\emph{\ semimodules} \cite{Alt1996}. One
reason for the interest in such notions is the phenomenon that assuming that 
\emph{all} semimodules over a given semiring $S$ are projective forces the
underlying semiring to be a (semisimple) \emph{ring} (cf., \cite[Theorem 3.4]%
{Ili2010}). Using a new notion of exact sequences of semimodules over a
semiring, Abuhlail \cite{Abu2014-CA} introduced the homological notion of 
\emph{exactly projective semimodules} ($e$-\emph{projective semimodules},
for short) assuming that an appropriate $Hom$functor preserves short exact
sequences (under the initial name of \emph{uniformly projective semimodules}%
).

The paper is divided into three sections.

\bigskip

In Section 1, we collect the basic definitions, examples and preliminaries
used in this paper. Among others, we include the definitions and basic
properties of \emph{exact sequences }introduced by Abuhlail \cite{Abu2014}.

\bigskip

In Section 2, we demonstrate the existence of \emph{pullbacks} (see \ref%
{pullb}) and \emph{pushouts} (Theorem \ref{existpushout}) in the category of
semimodules over an arbitrary semiring. Although no explicit construction of
the pushouts is given, we provide a description that is good enough to help
us in proving several results in the sequel.

\bigskip

In Section Three, we investigate mainly the $e$\emph{-projective semimodules}
over a semiring and clarify their relations with the notions of \emph{%
projective semimodules} as well as the so called $k$\emph{-projective
semimodules}. In Proposition \ref{proj->uproj}, we demonstrate that every
projective left semimodule is in fact $e$-projective. In Example \ref%
{e-proj-not-proj}, we show that the Boolean Algebra $\mathbb{B}$ considered
as a $\mathbb{Q}^{+}$-semimodule in the canonical way is $\mathbb{Q}^{+}$-$e$%
-projective but not $\mathbb{Q}^{+}$-projective. A complete characterization
of $k$-projective left semimodules through the right-splitting of short
exact sequences is given in Proposition \ref{char-k-proj}. In Lemma \ref%
{ret-proj} and Proposition \ref{dsum-e-proj}, we provide homological proofs
of the facts that the class of $e$-projective left $S$-semimodules is closed
under retracts and direct sums recovering part of \cite[Corollary 3.3]%
{AIKN2018}, where compact \emph{categorical proofs} were given.

\section{Preliminaries}

\label{prelim}

\qquad In this section, we provide the basic definitions and preliminaries
used in this work. Any notions that are not defined can be found in {our
main reference \cite{Gol1999}. We refer to \cite{Wis1991} for the
foundations of Module and Ring Theory.}

\begin{definition}
(\cite{Gol1999}) A \textbf{semiring}%
\index{Semiring} is a datum $(S,+,0,\cdot ,1)$ consisting of a commutative
monoid $(S,+,0)$ and a monoid $(S,\cdot ,1)$ such that $0\neq 1$ and%
\begin{eqnarray*}
a\cdot 0 &=&0=0\cdot a%
\text{ for all }a\in S; \\
a(b+c) &=&ab+ac\text{ and }(a+b)c=ac+bc\text{ for all }a,b,c\in S.
\end{eqnarray*}
\end{definition}

\begin{defns}
\label{def-semiring}(\cite{Gol1999}) Let $(S,+,0,\cdot ,1)$ be a semiring.

\begin{itemize}
\item If the monoid $(S,\cdot ,1)$ is commutative, we say that $S$ is a 
\emph{commutative semiring}.

\item The set of \emph{cancellative elements of }$S$ is defined as%
\begin{equation*}
K^{+}(S)=\{x\in S\mid x+y=x+z\Longrightarrow y=z\text{ for any }y,z\in S\}.
\end{equation*}%
We say that $S$ is a \emph{cancellative semiring} if $K^{+}(S)=S.$
\end{itemize}
\end{defns}

\begin{exs}
({\cite{Gol1999})}

\begin{itemize}
\item Every ring is a cancellative semiring.

\item Any \emph{distributive bounded lattice} $\mathcal{L}=(L,\vee ,1,\wedge
,0)$ is a commutative semiring.

\item Let $R$ be any ring. The set $\mathcal{I}=(Ideal(R),+,0\cdot ,R)$ of
ideals of $R$ is a semiring.

\item The sets $(\mathbb{Z}^{+},+,0,\cdot ,1)$ (resp. $(\mathbb{Q}%
^{+},+,0,\cdot ,1),$ $(\mathbb{Q}^{+},+,0,\cdot ,1)$) of non-negative
integers (resp. non-negative rational numbers, non-negative real numbers) is
a commutative cancellative semiring which is not a ring.

\item $M_{n}(S),$ the set of all $n\times n$ matrices over a semiring $S,$
is a semiring.

\item The Boolean algebra $\mathbb{B}:=\{0,1\}$ with $1+1=1$ is a semiring
called the \textbf{Boolean Semiring}.
\end{itemize}
\end{exs}

\begin{punto}
\cite{Gol1999} Let $S$ and $T$ be semirings. The categories $_{S}\mathbf{SM}$
of \textbf{left} $S$-\textbf{semimodules} with arrows the $S$-linear maps, $%
\mathbf{SM}_{T}$ of right $S$-semimodules with arrows the $T$-linear maps,
and $_{S}\mathbf{SM}_{T}$ of $(S,T)$-bisemimodules are defined in the usual
way (as for modules and bimodules over rings). We write $L\leq _{S}M$ to
mean that $M$ is a left (right) $S$-semimodule and $L$ is an $S$\textbf{%
-subsemimodule }of $M.$
\end{punto}

\begin{ex}
The category of $\mathbb{Z}^{+}$-semimodules is nothing but the category of
commutative monoids.
\end{ex}

\begin{ex}
Let $(S,+,0,\cdot ,1)$ be a semiring. Then $S$ and $S^{(\Lambda )}$ (the
direct sum of $S$ over a non-empty index set $\Lambda $) are $(S,S)$%
-bisemimodules with left and right actions induced by ``$\cdot $".
\end{ex}

\begin{ex}
(\cite[page 150, 154]{Gol1999}) Let $S$ be a semiring, $M$ be a left $S$%
-semimodule and $L\subseteq M.$ The \textbf{subtractive closure }of $L$ is
defined as%
\begin{equation}
\overline{L}:=\{m\in M\mid \text{ }m+l=l^{\prime }\text{ for some }%
l,l^{\prime }\in L\}.  \label{L-s-closure}
\end{equation}%
One can easily check that $\overline{L}=Ker(M\overset{\pi }{\longrightarrow }%
M/L),$ where $\pi $ is the canonical projection. We say that $L$ is \textbf{%
subtractive, }if $L=\overline{L}.$ The left $S$-semimodule $M$ is a \textbf{%
subtractive semimodule, }if every $S$-subsemimodule $L\leq _{S}M$ is
subtractive. 
If the only $S$-subsemimodules of $M$ are $\{0\}$ and $M,$ then we say that $%
M$ is \textbf{ideal-simple}.
\end{ex}

\begin{definition}
\cite[page 162]{Gol1999} Let $S$ be a semiring. An equivalence relation $%
\rho $ on a left $S$-semimodule $M$ is a \textbf{congruence relation,} if it
preserves the addition and the scalar multiplication on $M,$ \emph{i.e. }for
all $s\in S$ and $m,m^{\prime },n,n^{\prime }\in M:$%
\begin{equation*}
m\rho m^{\prime }\text{ and }n\rho n^{\prime }\Longrightarrow (m+m^{\prime
})\rho (n+n^{\prime }),
\end{equation*}%
\begin{equation*}
m\rho m^{\prime }\Longrightarrow (sm)\rho (sm^{\prime }).
\end{equation*}
\end{definition}

\begin{lem}
\label{s-char}A left $S$-semimodule $M$ is ideal-simple if and only if every
non-zero $S$-linear map to $M$ is surjective.
\end{lem}

\begin{punto}
\label{variety}(cf., \cite{AHS2004})\ The category $_{S}\mathbf{SM}$ of left
semimodules over a semiring $S$ is a \emph{variety} in the sense of
Universal Algebra (closed under homomorphic images, subobjects and arbitrary
products). Whence $_{S}\mathbf{SM}$ is complete, i.e. has all limits (e.g.,
direct products, equalizers, kernels, pullbacks, inverse limits) and
cocomplete, i.e. has all colimits (e.g., direct coproducts, coequalizers,
cokernels, pushouts, direct colimits).
\end{punto}

\begin{definition}
(\cite[page 184]{Gol1999}) Let $S$ be a semiring. A left $S$-semimodule $M$
is the \textbf{direct sum}%
\index{direct sum} of a family $\{L_{\lambda }\}_{\lambda \in \Lambda }$ of $%
S$-subsemimodules $L_{\lambda }\leq _{S}M,$ and we write $%
M=\bigoplus\limits_{\lambda \in \Lambda }L_{\lambda }$, if every $m\in M$
can be written in a \emph{unique way} as a finite sum $m=l_{\lambda
_{1}}+\cdots +l_{\lambda _{k}}$ where $l_{\lambda _{i}}\in L_{\lambda _{i}}$
for each $i=1,\cdots ,k.$ Equivalently, $M=\bigoplus\limits_{\lambda \in
\Lambda }L_{\lambda }$ if $M=\sum\limits_{\lambda \in \Lambda }L_{\lambda }$
and for each finite subset $A\subseteq \Lambda $ with $l_{a},l_{a}^{\prime
}\in L_{a},$ we have:%
\begin{equation*}
\sum\limits_{a\in A}l_{a}=\sum\limits_{a\in A}l_{a}^{\prime }\Longrightarrow
l_{a}=l_{a}^{\prime }%
\text{ for all }a\in A.
\end{equation*}
\end{definition}

\begin{punto}
An $S$-semimodule $N$ is a \textbf{retract} of an $S$-semimodule $M$ if
there exists a (surjective) $S$-linear map $\theta :M\longrightarrow N$ and
an (injective) $S$-linear map $\psi :N\longrightarrow M$ such that $\theta
\circ \psi =\mathrm{id}_{N}$ (equivalently, $N\simeq \alpha (M)$ for some 
\emph{idempotent} endomorphism $\alpha \in \mathrm{End}(M_{S})$).
\end{punto}

\begin{punto}
An $S$-semimodule $N$ is a \textbf{direct summand} of an $S$-semimodule $M$ (%
\emph{i.e.} $M=N\oplus N^{\prime }$ for some $S$-subsemimodule $N^{\prime }$
of $M$) if and only if there exists $\alpha \in \mathrm{Comp}(\mathrm{End}%
(M_{S}))$ s.t. $\alpha (M)=N$ where for any semiring $T$ we set%
\begin{equation*}
\mathrm{Comp}(T)=\{t\in T\mid \text{ }\exists \text{ }\widetilde{t}\in T%
\text{ with }t+\widetilde{t}=1_{T}\text{ and }t\widetilde{t}=0_{T}=%
\widetilde{t}t\}.
\end{equation*}%
Indeed, every direct summand of $M$ is a retract of $M;$ the converse is not
true in general; for example $N_{1}$ in Example \ref{nonseproj} is a retract
of $M_{2}(\mathbb{R^{+}})$ that is not a direct summand. Golan \cite[%
Proposition 16.6]{Gol1999} provided characterizations of direct summands.
\end{punto}

\begin{rems}
\label{d-iso}Let $M$ be a left $S$-semimodule and $K,L\leq _{S}M$ be $S$%
-semimodules of $M.$

\begin{enumerate}
\item If $K+L$ is direct, then $K\cap L=0.$ The converse is \emph{not} true
in general (for a counterexample, see Example see \ref{not-direct}).

\item If $M=K\oplus L$, then $M/K\simeq L.$
\end{enumerate}
\end{rems}

\begin{ex}
\label{not-direct}Let $S=M_{2}(\mathbb{R}^{+})$. Notice that%
\begin{equation*}
E_{1}=\left\{ \left[ {%
\begin{array}{cc}
a & 0 \\ 
b & 0%
\end{array}%
}\right] |\text{ }a,b\in \mathbb{R}^{+}\right\} \text{ and }N_{\geq
1}\left\{ \left[ {%
\begin{array}{cc}
a & c \\ 
b & d%
\end{array}%
}\right] |\text{ }a\leq c,b\leq d,a,b,c,d\in \mathbb{R}^{+}\right\}
\end{equation*}%
are left ideals of $S$ with $E_{1}\cap N_{\geq 1}=\{0\}$. However, the sum $%
E_{1}+N_{\geq 1}$ is not direct since 
\begin{equation*}
\left[ {%
\begin{array}{cc}
1 & 0 \\ 
0 & 0%
\end{array}%
}\right] +\left[ {%
\begin{array}{cc}
0 & 1 \\ 
0 & 0%
\end{array}%
}\right] =\left[ {%
\begin{array}{cc}
0 & 0 \\ 
0 & 0%
\end{array}%
}\right] +\left[ {%
\begin{array}{cc}
1 & 1 \\ 
0 & 0%
\end{array}%
}\right] .
\end{equation*}
\end{ex}

\subsection*{Exact Sequences}

\bigskip

Throughout, $(S,+,0,\cdot ,1)$ is a semiring and, unless otherwise
explicitly mentioned, an $S$-module is a \emph{left }$S$-semimodule.

\bigskip

\begin{definition}
A morphism of left $S$-semimodules $f:L\rightarrow M$ is

$k$-\textbf{normal}, if whenever $f(m)=f(m^{\prime })$ for some $m,m^{\prime
}\in M,$ we have $m+k=m^{\prime }+k^{\prime }$ for some $k,k^{\prime }\in
Ker(f);$

$i$-\textbf{normal}, if $\func{Im}(f)=\overline{f(L)}$ ($:=\{m\in M|\text{ }%
m+l\in L\text{ for some }l\in L\}$).

\textbf{normal}, if $f$ is both $k$-normal and $i$-normal.
\end{definition}

\begin{rems}
\begin{enumerate}
\item Among others, Takahashi (\cite{Tak1981}) and Golan \cite{Gol1999}
called $k$-normal (resp., $i$-normal, normal) $S$-linear maps $k$\emph{%
-regular} (resp., $i$\emph{-regular}, \emph{regular}) morphisms. We changed
the terminology to avoid confusion with the regular monomorphisms and
regular epimorphisms in Category Theory which have different meanings when
applied to categories of semimodules.

\item Our terminology is consistent with Category Theory noting that: every
surjective $S$-linear map is $i$-normal, whence the $k$-normal surjective $S$%
-linear map are normal and are precisely the so-called \textbf{normal
epimorphisms}. On the other hand, the injective $S$-linear maps are $k$%
-normal, whence the $i$-normal injective $S$-linear maps are normal and are
precisely the so called \textbf{normal monomorphisms} (see \cite{Abu2014}).
\end{enumerate}
\end{rems}

The following technical lemma is helpful in several proofs in this and
forthcoming related papers.

\begin{lem}
\label{i-normal}Let $L\overset{f}{\rightarrow }M\overset{g}{\rightarrow }N$
be a sequence of semimodules.

\begin{enumerate}
\item Let $g$ be injective.

\begin{enumerate}
\item $f$ is $k$-normal if and only if $g\circ f$ is $k$-normal.

\item If $g\circ f$ is $i$-normal (normal), then $f$ is $i$-normal (normal).

\item Assume that $g$ is $i$-normal. Then $f$ is $i$-normal (normal) if and
only if $g\circ f$ is $i$-normal (normal).
\end{enumerate}

\item Let $f$ be surjective.

\begin{enumerate}
\item $g$ is $i$-normal if and only if $g\circ f$ is $i$-normal.

\item If $g\circ f$ is $k$-normal (normal), then $g$ is $k$-normal (normal).

\item Assume that $f$ is $k$-normal. Then $g$ is $k$-normal (normal) if and
only if $g\circ f$ is $k$-normal (normal).
\end{enumerate}
\end{enumerate}
\end{lem}

\begin{Beweis}
\begin{enumerate}
\item Let $g$ be injective; in particular, $g$ is $k$-normal.

\begin{enumerate}
\item Assume that $f$ is $k$-normal. Suppose that $(g\circ f)(l_{1})=(g\circ
f)(l_{2})$ for some $l_{1},l_{2}\in L.$ Since $g$ is injective, $%
f(l_{1})=f(l_{2}).$ By assumption, there exist $k_{1},k_{2}\in \mathrm{Ker}%
(f)$ such that $l_{1}+k_{1}=l_{2}+k_{2}.$ Since $\mathrm{Ker}(f)\subseteq 
\mathrm{Ker}(g\circ f),$ we conclude that $g\circ f$ is $k$-normal. On the
other hand, assume that $g\circ f$ is $k$-normal. Suppose that $%
f(l_{1})=f(l_{2})$ for some $l_{1},l_{2}\in L.$ Then $(g\circ
f)(l_{1})=(g\circ f)(l_{2})$ and so there exist $k_{1},k_{2}\in \mathrm{Ker}%
(g\circ f)$ such that $l_{1}+k_{1}=l_{2}+k_{2}.$ Since $g$ is injective, $%
\mathrm{Ker}(g\circ f)=\mathrm{Ker}(f)$ whence $f$ is $k$-normal.

\item Assume that $g\circ f$ is $i$-normal. Let $m\in \overline{f(L)},$ so
that $m+f(l_{1})=f(l_{2})$ for some $l_{1},l_{2}\in L.$ Then $g(m)\in 
\overline{(g\circ f)(L)}=(g\circ f)(L).$ Since $g$ is injective, $m\in f(L).$
So, $f$ is $i$-normal.

\item Assume that $g$ and $f$ are $i$-normal. Let $n\in \overline{(g\circ
f)(L)},$ so that $n+g(f(l_{1}))=g(f(l_{2}))$ for some $l_{1},l_{2}\in L.$
Since $g$ is $i$-normal, $n\in g(M)$ say $n=g(m)$ for some $m\in M.$ But $g$
is injective, whence $m+f(l_{1})=f(l_{2}),$ i.e. $m\in \overline{f(L)}=f(L)$
since $f$ is $i$-normal by assumption. So, $n=g(m)\in (g\circ f)(L).$ We
conclude that $g\circ f$ is $i$-normal.
\end{enumerate}

\item Let $f$ be surjective; in particular, $f$ is $i$-normal.

\begin{enumerate}
\item Assume that $g$ is $i$-normal. Let $n\in \overline{(g\circ f)(L)}$ so
that $n+g(f(l_{1}))=g(f(l_{2}))$ for some $l_{1},l_{2}\in L.$ Since $g$ is $%
i $-normal, $n=g(m)$ for some $m\in M.$ Since $f$ is surjective, $n=g(m)\in
(g\circ f)(L).$ So, $g\circ f$ is $i$-normal.

On the other hand, assume that $g\circ f$ is $i$-normal. Let $n\in \overline{%
g(M)},$ so that $n+g(m_{1})=g(m_{2})$ for some $m_{1},m_{2}\in M.$ Sine $f$
is surjective, there exist $l_{1},l_{2}\in L$ such that $f(l_{1})=m_{1}$ and 
$f(l_{2})=m_{2}.$ Then, $n+(g\circ f)(l_{1})=(g\circ f)(l_{2}),$ i.e. $n\in 
\overline{(g\circ f)(L)}=(g\circ f)(L)\subseteq g(M).$ So, $g$ is $i$-normal.

\item Assume that $g\circ f$ is $k$-normal. Suppose that $g(m_{1})=g(m_{2})$
for some $m_{1},m_{2}\in M.$ Since $f$ is surjective, we have $(g\circ
f)(l_{1})=(g\circ f)(l_{2})$ for some $l_{1},l_{2}\in L.$ By assumption, $%
g\circ f$ is $k$-normal and so there exist $k_{1},k_{2}\in \mathrm{Ker}%
(g\circ f)$ such that $l_{1}+k_{1}=l_{2}+k_{2}$ whence $%
m_{1}+f(k_{1})=m_{2}+f(k_{2}).$ Indeed, $f(k_{1}),f(k_{2})\in \mathrm{Ker}%
(g),$ i.e. $g$ is $k$-normal.

\item Assume that $f$ and $g$ are $k$-normal. Suppose that $(g\circ
f)(l_{1})=(g\circ f)(l_{2})$ for some $l_{1},l_{2}\in L.$ Since $g$ is $k$%
-normal, we have $f(l_{1})+k_{1}=f(l_{2})+k_{2}$ for some $k_{1},k_{2}\in 
\mathrm{Ker}(g).$ But $f$ is surjective; whence $k_{1}=f(l_{1}^{\prime })$
and $k_{2}=f(l_{2}^{\prime })$ for some $l_{1}^{\prime },l_{2}^{\prime }\in
L,$ i.e. $f(l_{1}+l_{1}^{\prime })=f(l_{2}+l_{2}^{\prime }).$ Since $f$ is $%
k $-normal, $l_{1}+l_{1}^{\prime }+k_{1}^{\prime }=l_{2}+l_{2}^{\prime
}+k_{2}^{\prime }$ for some $k_{1}^{\prime },k_{2}^{\prime }\in \mathrm{Ker}%
(f).$ Indeed, $l_{1}^{\prime }+k_{1}^{\prime },$ $l_{2}^{\prime
}+k_{2}^{\prime }\in \mathrm{Ker}(g\circ f).$ We conclude that $g\circ f$ is 
$k$-normal.$\blacksquare $
\end{enumerate}
\end{enumerate}
\end{Beweis}

There are several notions of exactness for sequences of semimodules. In this
paper, we use the relatively new notion introduced by Abuhlail:

\begin{definition}
\label{Abu-exs}(\cite[2.4]{Abu2014}) A sequence 
\begin{equation}
L\overset{f}{\longrightarrow }M\overset{g}{\longrightarrow }N  \label{LMN}
\end{equation}%
of left $S$-semimodules is \textbf{exact}, if $g$ is $k$-normal and $%
f(L)=Ker(g).$
\end{definition}

\begin{punto}
\label{def-exact}We call a sequence of $S$-semimodules $L\overset{f}{%
\rightarrow }M\overset{g}{\rightarrow }N$

\emph{proper-exact} if $f(L)=\mathrm{Ker}(g)$ (exact in the sense of
Patchkoria \cite{Pat2003});

\emph{semi-exact} if $\overline{f(L)}=\mathrm{Ker}(g)$ (exact in the sense
of Takahashi \cite{Tak1981});

\emph{quasi-exact} if $\overline{f(L)}=\mathrm{Ker}(g)$ and $g$ is $k$%
-normal (exact in the sense of Patil and Doere \cite{PD2006}).
\end{punto}

\begin{punto}
We call a (possibly infinite) sequence of $S$-semimodules 
\begin{equation}
\cdots \rightarrow M_{i-1}\overset{f_{i-1}}{\rightarrow }M_{i}\overset{f_{i}}%
{\rightarrow }M_{i+1}\overset{f_{i+1}}{\rightarrow }M_{i+2}\rightarrow \cdots
\label{chain}
\end{equation}

\emph{chain complex} if $f_{j+1}\circ f_{j}=0$ for every $j;$

\emph{exact} (resp., \emph{proper-exact}, \emph{semi-exact, quasi-exact}) if
each partial sequence with three terms $M_{j}\overset{f_{j}}{\rightarrow }%
M_{j+1}\overset{f_{j+1}}{\rightarrow }M_{j+2}$ is exact (resp.,
proper-exact, semi-exact, quasi-exact).

A \textbf{short exact sequence}%
\index{short exact sequence} (or a \textbf{Takahashi extension}%
\index{Takahashi extension} \cite{Tak1982b}) of $S$-semimodules is an exact
sequence of the form%
\begin{equation*}
0\longrightarrow L\overset{f}{\longrightarrow }M\overset{g}{\longrightarrow }%
N\longrightarrow 0
\end{equation*}
\end{punto}

\begin{rem}
In the sequence (\ref{LMN}), the inclusion $f(L)\subseteq Ker(g)$ forces $%
f(L)\subseteq 
\overline{f(L)}\subseteq Ker(g),$ whence the assumption $f(L)=Ker(g)$
guarantees that $f(L)=\overline{f(L)},$ \emph{i.e.} $f$ is $i$-normal. So,
the definition puts conditions on $f$ and $g$ that are dual to each other
(in some sense).
\end{rem}

The following result shows some of the advantages of the Abuhlail's
definition of exact sequences over the previous ones:

\begin{lem}
\label{exact}Let $L,M$ and $N$ be $S$-semimodules.

\begin{enumerate}
\item $0\longrightarrow L\overset{f}{\longrightarrow }M$ is exact if and
only if $f$ is injective.

\item $M\overset{g}{\longrightarrow }N\longrightarrow 0$ is exact if and
only if $g$ is surjective.

\item $0\longrightarrow L\overset{f}{\longrightarrow }M\overset{g}{%
\longrightarrow }N$ is semi-exact and $f$ is normal (proper-exact and $f$ is
normal) if and only if $L\simeq \mathrm{Ker}(g).$

\item $0\longrightarrow L\overset{f}{\longrightarrow }M\overset{g}{%
\longrightarrow }N$ is exact if and only if $L\simeq \mathrm{Ker}(g)$ and $g$
is $k$-normal.

\item $L\overset{f}{\longrightarrow }M\overset{g}{\longrightarrow }%
N\longrightarrow 0$ is semi-exact and $g$ is normal if and only if $N\simeq
M/f(L).$

\item $L\overset{f}{\longrightarrow }M\overset{g}{\longrightarrow }%
N\longrightarrow 0$ is exact if and only if $N\simeq M/f(L)$ and $f$ is $i$%
-normal.

\item $0\longrightarrow L\overset{f}{\longrightarrow }M\overset{g}{%
\longrightarrow }N\longrightarrow 0$ is exact if and only if $L\simeq 
\mathrm{Ker}(g)$ and $N\simeq M/L.$
\end{enumerate}
\end{lem}

\begin{cor}
\label{M/L}The following assertions are equivalent:

\begin{enumerate}
\item $0\rightarrow L\overset{f}{\rightarrow }M\overset{g}{\rightarrow }%
N\rightarrow 0$ is an exact sequence of $S$-semimodules;

\item $L\simeq \mathrm{Ker}(g)$ and $N\simeq M/f(L)$;

\item $f$ is injective, $f(L)=\mathrm{Ker}(g),$ $g$ is surjective and ($k$%
-)normal.

In this case, $f$ and $g$ are normal morphisms.
\end{enumerate}
\end{cor}

\begin{rem}
A morphism of semimodules $\gamma :X\longrightarrow Y$ is an isomorphism if
and only if $0\longrightarrow X\overset{\gamma }{\longrightarrow }%
Y\longrightarrow 0$ is exact if and only if $\gamma $ is a normal bimorphism
(\emph{i.e.} $\gamma $ is a normal monomorphism and a normal epimorphism).
The assumption on $\gamma $ to be normal cannot be removed here. For
example, the embedding $\iota :\mathbb{Z}^{+}\longrightarrow \mathbb{Z}$ is
a bimorphism of commutative monoids ($\mathbb{Z}^{+}$-semimodules) which is
not an isomorphism. Notice that $\iota $ is not $i$-normal; in fact $%
\overline{\iota (\mathbb{Z}^{+}})=\mathbb{Z}.$
\end{rem}

\begin{rem}
An $S$-linear map is a monomorphism if and only if it is injective. Every
surjective $S$-linear map is an epimorphism. The converse is not true in
general.
\end{rem}

\begin{ex}
The embedding $\iota:\mathbb{Z}^+\rightarrow \mathbb{Z}$ is a monoid
epimorphism as $(f \circ \iota)(1_{\mathbb{Z}^+})=(g\circ\iota)(1_{\mathbb{Z}%
^+})$ implies $f(1_{\mathbb{Z}})=g(1_{\mathbb{Z}})$ and $f=g$ for every
monoid morphisms $f,g:\mathbb{Z}\rightarrow M$. However, it is clear that $%
\iota$ is not surjective.
\end{ex}

\begin{lem}
\label{1st-IT}\emph{(Compare with \cite[Proposition 4.3.]{Tak1981})} Let $%
\gamma :X\rightarrow Y$ be a morphism of $S$-semimodules.

\begin{enumerate}
\item The sequence%
\begin{equation}
0\rightarrow \mathrm{Ker}(\gamma )\overset{\mathrm{\ker }(\gamma )}{%
\longrightarrow }X\overset{\gamma }{\rightarrow }Y\overset{\mathrm{\mathrm{%
coker}}(\gamma )}{\longrightarrow }\mathrm{\mathrm{\mathrm{\mathrm{\mathrm{%
Coker}}}}}(\gamma )\rightarrow 0  \label{ker-coker}
\end{equation}

with canonical $S$-linear maps is semi-exact. Moreover, (\ref{ker-coker}) is
exact if and only if $\gamma $ is normal.

\item We have two exact sequences%
\begin{equation*}
0\rightarrow \overline{\gamma (X)}\overset{\mathrm{ker}(\mathrm{\mathrm{coker%
}}(\gamma ))}{\longrightarrow }Y\overset{\mathrm{\mathrm{coker}}(\gamma )}{%
\longrightarrow }Y/\gamma (X)\rightarrow 0.
\end{equation*}%
and%
\begin{equation*}
0\rightarrow \mathrm{Ker}(\gamma )\overset{\mathrm{ker}(\gamma )}{%
\longrightarrow }X\overset{\mathrm{\mathrm{coker}}(\mathrm{ker}(\gamma ))}{%
\longrightarrow }X/\mathrm{Ker}(\gamma )\rightarrow 0.
\end{equation*}
\end{enumerate}
\end{lem}

\begin{cor}
\label{reg-sub}\emph{(Compare with \cite[Proposition 4.8.]{Tak1981}) }Let $M$
be an $S$-semimodule.

\begin{enumerate}
\item Let $\rho $ an $S$-congruence relation on $M$ and consider the
sequence of $S$-semimodules%
\begin{equation*}
0\longrightarrow \mathrm{Ker}(\pi _{\rho })\overset{\iota _{\rho }}{%
\longrightarrow }M\overset{\rho }{\longrightarrow }M/\rho \longrightarrow 0.
\end{equation*}

\begin{enumerate}
\item $0\rightarrow \mathrm{Ker}(\pi _{\rho })\overset{\iota _{\rho }}{%
\longrightarrow }M\overset{\pi _{\rho }}{\longrightarrow }M/\rho \rightarrow
0$ is exact.

\item $M/\rho =\mathrm{Coker}(\iota _{\rho })$.
\end{enumerate}

\item Let $L$ be an $S$-subsemimodule of $M$.

\begin{enumerate}
\item The sequence $0\rightarrow L\overset{\iota }{\longrightarrow }M\overset%
{\pi _{L}}{\longrightarrow }M/L\rightarrow 0$ is semi-exact.

\item $0\rightarrow \overline{L}\overset{\iota }{\longrightarrow }M\overset{%
\pi _{L}}{\longrightarrow }M/L\rightarrow 0$ is exact.

\item The following assertions are equivalent:

\begin{enumerate}
\item $0\rightarrow L\overset{\iota }{\longrightarrow }M\overset{\pi _{L}}{%
\longrightarrow }M/L\rightarrow 0$ is exact;

\item $L= \mathrm{Ker}(\pi _{L});$

\item $0\longrightarrow L\overset{\iota }{\longrightarrow }\overline{L}%
\longrightarrow 0$ is exact;

\item $L$ is a subtractive subsemimodule.
\end{enumerate}
\end{enumerate}
\end{enumerate}
\end{cor}

\begin{prop}
\label{adj-lim}\emph{(}cf., \emph{\cite[Proposition 3.2.2]{Bor1994})} Let $%
\mathfrak{C},\mathfrak{D}$ be arbitrary categories and $\mathfrak{C}\overset{%
F}{\longrightarrow }\mathfrak{D}\overset{G}{\longrightarrow }\mathfrak{C}$
be functors such that $(F,G)$ is an adjoint pair.

\begin{enumerate}
\item $F$ preserves all colimits which turn out to exist in $\mathfrak{C}.$

\item $G$ preserves all limits which turn out to exist in $\mathfrak{D}.$
\end{enumerate}
\end{prop}

\begin{cor}
\label{ad-l-cor}Let $S,$ $T$ be semirings and $_{T}F_{S}$ a $(T,S)$%
-bisemimodule. The covariant functor $\mathrm{Hom}_{T}(F,-):$ $_{T}\mathbf{SM%
}\longrightarrow $ $_{S}\mathbf{SM}$ preserves all limits.

\begin{enumerate}
\item For every family of left $T$-semimodules $\{Y_{\lambda }\}_{\Lambda },$
we have a canonical isomorphism of left $S$-semimodules%
\begin{equation*}
\mathrm{Hom}_{T}(F,\prod\limits_{\lambda \in \Lambda }Y_{\lambda })\simeq
\prod\limits_{\lambda \in \Lambda }\mathrm{Hom}_{T}(F,Y_{\lambda }).
\end{equation*}

\item For any inverse system of left $T$-semimodules $(X_{j},\{f_{jj^{\prime
}}\})_{J},$ we have an isomorphism of left $S$-semimodules%
\begin{equation*}
\mathrm{Hom}_{T}(F,\lim_{\longleftarrow }X_{j})\simeq \text{ }%
\lim_{\longleftarrow }\mathrm{Hom}_{T}(F,X_{j}).
\end{equation*}

\item $\mathrm{Hom}_{T}(F,-)$ preserves equalizers;

\item $\mathrm{Hom}_{T}(F,-)$ preserves kernels.
\end{enumerate}
\end{cor}

\begin{Beweis}
The proof can be obtained as a direct consequence of Proposition \ref%
{adj-lim} and the fact that $(F\otimes _{S}-,\mathrm{Hom}_{T}(F,-))$ is an
adjoint pair of covariant functors \cite{KN2011}.$\blacksquare $
\end{Beweis}

Corollary \ref{ad-l-cor} allows us to improve \cite[Theorem 2.6]{Tak1982a}.

\begin{prop}
\label{lr-exact}Let $_{T}G_{S}$ be $(T,S)$-bisemimodule and consider the
functor $\mathrm{Hom}_{T}(G,-):$ $_{T}\mathbf{SM}\longrightarrow $ $_{S}%
\mathbf{SM}.$ Let%
\begin{equation}
0\longrightarrow L\overset{f}{\rightarrow }M\overset{g}{\rightarrow }N
\label{lr}
\end{equation}%
be a sequence of left $T$-semimodules and consider the following sequence of
left $S$-semimodules%
\begin{equation}
0\longrightarrow \mathrm{Hom}_{T}(G,L)\overset{(G,f)}{\rightarrow }\mathrm{%
Hom}_{T}(G,M)\overset{(G,g)}{\longrightarrow }\mathrm{Hom}_{T}(G,N).
\label{GL}
\end{equation}

\begin{enumerate}
\item If the sequence $0\longrightarrow L\overset{f}{\rightarrow }M$ is
exact and $f$ is normal, then 
\begin{equation*}
0\longrightarrow \mathrm{Hom}_{T}(G,L)\overset{(G,f)}{\rightarrow }\mathrm{%
Hom}_{T}(G,M)
\end{equation*}%
is exact and $(G,f)$ is normal.

\item If (\ref{lr})\emph{\ }is proper-exact and $f$ is normal, then (\ref{GL}%
)\emph{\ }is proper exact and $(G,f)$ is normal.

\item If (\ref{lr})\emph{\ }is exact and $\mathrm{Hom}_{T}(G,-)$ preserves $%
k $-normal morphisms, then (\ref{GL})\emph{\ }is exact.
\end{enumerate}
\end{prop}

\begin{Beweis}
\begin{enumerate}
\item The following implications are obvious: $0\longrightarrow L\overset{f}{%
\rightarrow }M$ is exact $\Longrightarrow $ $f$ is injective $%
\Longrightarrow $ $(G,f)$ is injective $\Longrightarrow 0\longrightarrow 
\mathrm{Hom}_{T}(G,L)\overset{(G,f)}{\rightarrow }\mathrm{Hom}_{T}(G,M)$ is
exact. Assume that $f$ is normal and consider the short exact sequence of $S$%
-semimodules%
\begin{equation*}
0\longrightarrow L\overset{f}{\longrightarrow }M\overset{\pi _{L}}{%
\longrightarrow }M/L\longrightarrow 0.
\end{equation*}%
Notice that $L=\mathrm{Ker}(\pi _{L})$ by Lemma \ref{exact}. By Corollary %
\ref{ad-l-cor}, $\mathrm{Hom}_{T}(G,-)$ preserves kernels and so $(G,f)=%
\mathrm{ker}(G,\pi _{L})$ whence normal.

\item Apply Lemma \ref{exact} (3): The proper-exactness of (\emph{\ref{lr}})
and the normality of $f$ are equivalent to $L\simeq \mathrm{Ker}(g).$ Since $%
\mathrm{Hom}_{T}(G,-)$ preserves kernels, we deduce that $\mathrm{Hom}%
_{T}(G,L)=\mathrm{Ker}((G,g))$, whence (\ref{GL})\emph{\ }is proper-exact
and $(G,f)$ is normal by Lemma \ref{exact} (3).

\item The statement follows directly from (2) and the assumption on $\mathrm{%
Hom}_{T}(G,-).\blacksquare $
\end{enumerate}
\end{Beweis}

\begin{punto}
\label{gam-adj}Let $\gamma :T\longrightarrow S$ be a morphism of semirings.
Then we have an \emph{adjoint pair of functors} $(F(X),Hom_{T}(S,-)),$ where 
$F(X)=X$ with the action $tx=\gamma (t)x$ for all $t\in T$ and $x\in X$ and $%
(s_{1}f)(s)=f(ss_{1})$ for all $s_{1},s\in S$ and $f\in Hom_{T}(S,Y)$ for
every left $T$-semimodule $Y.$ In particular, we have for all $X\in $ $_{S}%
\mathbf{SM}$ and $Y\in $ $_{T}\mathbf{SM}$ a natural \emph{isomorphism} of
commutative monoids%
\begin{equation}
\theta _{X,Y}:Hom_{S}(X,Hom_{T}(S_{S},Y))\longrightarrow Hom_{T}(X,Y),\text{ 
}f\mapsto \lbrack x\mapsto f(x)(1_{S})]  \label{tht}
\end{equation}%
with inverse%
\begin{equation}
\phi _{X,Y}:Hom_{T}(X,Y)\longrightarrow Hom_{S}(X,Hom_{T}(S_{S},Y)),\text{ }%
g\mapsto \lbrack x\mapsto [s\mapsto g(sx)]].  \label{tht-1}
\end{equation}
\end{punto}

\section{Pullbacks and Pushouts}

Throughout, $(S,+,0,\cdot ,1)$ is a semiring and, unless otherwise
explicitly mentioned, an $S$-module is a \emph{left }$S$-semimodule. The
category of left $S$-semimodules is denoted by $_{S}\mathbf{SM}$.

\bigskip

The category $_{S}\mathbf{SM}$ of left $S$-semimodules has pullbacks and
pushouts.

\bigskip

The pullbacks in $_{S}\mathbf{SM}$ are constructed in a way similar to that
of pullbacks in the category of modules over a ring.

\begin{punto}
\label{pullb}(\cite[1.7]{Tak1982b}) Let $f:A\rightarrow C$ and $%
g:B\rightarrow C$ be morphisms of left $S$-semimodules. The \textbf{pullback 
}%
\index{pullback}of $(f,g)$ is $(Q;f^{\prime },g^{\prime }),$ where%
\begin{eqnarray}
Q &:&=\{(a,b)\in A\times B|%
\text{ }f(a)=g(b)\}  \label{Pull} \\
g^{\prime } &:&Q\rightarrow A,\text{ }(a,b)\mapsto a;  \notag \\
f^{\prime } &:&Q\rightarrow C,\text{ }(a,b)\mapsto b,  \notag
\end{eqnarray}%
\begin{equation}
\xymatrix{Q^\ast \ar@{-->}[rd]^{\varphi} \ar@/^-1.0pc/[rdd]_{f^\ast}
\ar@/^1.0pc/[rrd]^{g^\ast} & & \\ & Q \ar[r]^{g^\prime} \ar[d]_{f^\prime} &
A \ar[d]^{f} \\ & B \ar[r]_{g} & C }  \label{pbc}
\end{equation}%
and whenever $(Q^{\ast };f^{\ast },g^{\ast })$ satisfies $f^{\ast }\circ
g=g^{\ast }\circ f$, there exists a unique $S$-linear map $\varphi :Q^{\ast
}\rightarrow Q$ such that $f\circ \varphi =f^{\ast }$ and $g\circ \varphi
=g^{\ast }$.
\end{punto}

Although the existence of pushouts in the category $_{S}\mathbf{SM}$ is
guaranteed since this category is a \emph{variety} in the sense of Universal
Algebra (see \ref{variety}), the construction of pushouts in it is much more
subtle that the construction of pushouts in the category of modules over a
ring (mainly because of the lack of subtraction).

\bigskip

This made some authors consider a special version of pushouts, e.g.,
Takahashi \cite{Tak1982b} who constructed in the so called $C$\textbf{%
-pushouts}, which coincide with the pushouts in the subcategory of \emph{%
cancellative} semimodules.

\begin{punto}
\label{cpush}(\cite[1.8]{Tak1982b}) Let $f:L\rightarrow M$ and $%
g:L\rightarrow N$ be morphisms of left $S$-semimodules. Consider the
congruence $\sim $ on $M\oplus N$ defined as 
\begin{equation}
(m_{1},n_{1})\sim (m_{2},n_{2})\Leftrightarrow \text{ }\exists \text{ }%
l_{1},l_{2}\in L:\text{ }m_{1}+f(l_{1})=m_{2}+f(l_{2})\text{ and }%
n_{1}+g(l_{2})=n_{2}+g(l_{1}).  \label{CP-cong}
\end{equation}%
The $C$\emph{-pushout} of $(f,g)$ is%
\begin{eqnarray}
CP &:&=(\iota _{M},\iota _{N};(M\oplus N)/\sim );  \label{CP} \\
\iota _{M} &:&M\rightarrow CP,\text{ }m\mapsto \lbrack (m,0)];  \notag \\
\iota _{N} &:&N\rightarrow CP,\text{ }n\mapsto \lbrack (0,n)].  \notag
\end{eqnarray}
\end{punto}

While the $C$-pushouts coincide with the natural pushout in the subcategory $%
_{S}\mathbf{CSM}$ of cancellative left semimodules, they fail to have the 
\emph{universal property} of pushouts in $_{S}\mathbf{SM}$.

\bigskip

In what follows, we demonstrate the construction of pushouts in $S$%
-semimodules $_{S}\mathbf{SM}$. The \emph{constructive proof} is the
objective of the following theorem which is already known to be true.

\begin{thm}
\label{existpushout}Let $f:L\rightarrow M$ and $g:L\rightarrow N$ be
morphisms of left $S$-semimodules. Then $(f,g)$ has a pushout.%

\begin{Beweis}
Consider%
\begin{equation*}
\begin{array}{cc}
\mathcal{P}:= & \{(g^{\prime },f^{\prime },P)\text{ }|\text{ }P\in \text{ }%
_{S}\mathbf{SM},\text{ }g^{\prime }:M\rightarrow P,f^{\prime }:N\rightarrow
P,\text{ }g^{\prime }\circ f=f^{\prime }\circ g, \\ 
& \pi _{(g^{\prime },f^{\prime })}:M\oplus N\longrightarrow P,\text{ }%
(m,n)\mapsto g^{\prime }(m)+f^{\prime }(n)\text{ is surjective}\}.%
\end{array}%
\end{equation*}%
\begin{equation*}
\xymatrix{L \ar[r]^{f} \ar[d]_{g} & M \ar[d]^{g'}\\ N \ar[r]_{f'} & P &
M\oplus N \ar[l]^{\pi_{(g',f')}}}
\end{equation*}%
Notice that $\mathcal{P}$ is not empty as $(0,0,0)\in \mathcal{P}$.

Define a relation $\leq $ on $\mathcal{P}$ as $(\tilde{g},\tilde{f},U)\leq
(f^{\prime },g^{\prime },P)$ if there exists an $S\text{-linear map }\alpha
:P\rightarrow U\text{ such that }\alpha \circ \pi _{(g^{\prime },f^{\prime
})}=\pi _{(\tilde{g},\tilde{f})}$, i.e. the following diagram is commutative%
\begin{equation*}
{\xymatrix{M\oplus N \ar[r]^{\pi_{(g',f')}}
\ar[d]_{\pi_{(\tilde{g},\tilde{f})}} & P \ar@{-->}[ld]^{\alpha}\\ U }}
\end{equation*}

\textbf{Step I:}\ $\mathcal{P}$ has a largest element $(\pi _{M},\pi _{N},%
\mathbf{P}),$ where%
\begin{eqnarray*}
\mathbf{P} &:&=(M\oplus N)/\rho , \\
(m_{1},n_{1})\rho (m_{2},n_{2}) &\Leftrightarrow &g_{\lambda
}(m_{1})+f_{\lambda }(n_{1})=g_{\lambda }(m_{2})+f_{\lambda }(n_{2})\text{ }%
\forall (g_{\lambda },f_{\lambda },P_{\lambda })\in \mathcal{P}; \\
\pi _{M} &:&M\longrightarrow (M\oplus N)/\rho ,\text{ }m\mapsto \lbrack
(m,0)]; \\
\pi _{N} &:&M\longrightarrow (M\oplus N)/\rho ,\text{ }n\mapsto \lbrack
(0,n)].
\end{eqnarray*}

\begin{itemize}
\item Notice that $(\pi _{M},\pi _{N},\mathbf{P})\in \mathcal{P}:$ for any $%
l\in L,$ we have for any $(g_{\lambda },f_{\lambda },P_{\lambda })\in 
\mathcal{P}$:%
\begin{equation*}
(g_{\lambda }\circ f)(l)+f_{\lambda }(0_{N})=(g_{\lambda }\circ
f)(l)=(f_{\lambda }\circ g)(l)=g_{\lambda }(0_{M})+(f_{\lambda }\circ g)(l)
\end{equation*}%
whence (by the definition of $\rho $):%
\begin{equation*}
(\pi _{M}\circ f)(l)=[(f(l),0)]_{\rho }=[(0,g(l))]_{\rho }=(\pi _{N}\circ
g)(l).
\end{equation*}

\item For every $(g_{\lambda },f_{\lambda },P_{\lambda })\in \mathcal{P},$
consider the $S$-linear map%
\begin{equation*}
\alpha _{\lambda }:\mathbf{P}\longrightarrow P_{\lambda },\text{ }%
[(m,n)]_{\rho }\mapsto g_{\lambda }(m)+f_{\lambda }(n)
\end{equation*}%
Notice that $\alpha _{\lambda }$ is well defined: if $[(m_{1},n_{1})]_{\rho
}=[(m_{2},n_{2})]_{\rho },$ then it follows by the definition of $\rho $
that 
\begin{equation*}
\alpha ([(m_{1},n_{1})]_{\rho })=g_{\lambda }(m_{1})+f_{\lambda
}(n_{1})=g_{\lambda }(m_{2})+f_{\lambda }(n_{2})=\alpha
([(m_{2},n_{2})]_{\rho }).
\end{equation*}%
Moreover, the following diagram 
\begin{equation*}
\xymatrix{M\oplus N \ar[r]^{\pi_{(\pi_M,\pi_N)}}
\ar[d]_{\pi_{(g_\lambda,f_\lambda)}} & {\bf P}
\ar@{-->}[ld]^{\alpha_\lambda}\\ P_\lambda }
\end{equation*}%
is commutative:\ indeed, for all $(m,n)\in M\oplus N$ we have%
\begin{equation*}
(\alpha _{\lambda }\circ \pi _{(\pi _{M},\pi _{N})})((m,n))=\alpha _{\lambda
}[(m,n)]_{\rho }=g_{\lambda }(m)+f_{\lambda }(n)=\pi _{(g_{\lambda
},f_{\lambda })}(m,n).
\end{equation*}
\end{itemize}

\textbf{Step II:}\ A largest element $(g^{\prime },f^{\prime };P)$ of $%
\mathcal{P}$ is a pushout of $(f,g).$ By the definition of $\mathcal{P}$, we
have $g^{\prime }\circ f=f^{\prime }\circ g.$ So it remains to prove the it
has the universal property of pushouts.

\begin{itemize}
\item Let $Q$ be a left $S$-semimodule along with $S$-linear maps $g^{\ast
}:M\rightarrow Q$ and $f^{\ast }:N\rightarrow Q$ satisfying $g^{\ast }\circ
f=f^{\ast }\circ g$. Since $\pi _{(g^{\prime },f^{\prime })}$ is surjective,
there exists for each $p\in P$ some $(m,n)\in M\oplus N,$ such that $%
p=g^{\prime }(m)+f^{\prime }(n)$. Define%
\begin{equation*}
\varphi :P\rightarrow Q,\text{ }p\mapsto g^{\ast }(m)+f^{\ast }(n).
\end{equation*}%
\begin{equation}
\xymatrix{L \ar[r]^{f} \ar[d]_{g} & M \ar[d]^{g'} \ar@/^1.0pc/[rdd]^{g^*}\\
N \ar[r]_{f'} \ar@/^-1.0pc/[rrd]^{f^*} & P \ar@{-->}[rd]^{\varphi} \\ & & Q}
\label{UP}
\end{equation}

\item It follows directly from the definition of $\varphi $ that $\varphi
\circ g^{\prime }=g^{\ast }$ and $\varphi \circ f^{\prime }=f^{\ast }.$

\textbf{Claim:} $\varphi $ is \emph{well defined}. Suppose that there exist $%
(m_{1},n_{1}),$ $(m_{2},n_{2})\in M\oplus N$ such that $g^{\prime
}(m_{1})+f^{\prime }(n_{1})=p=g^{\prime }(m_{2})+f^{\prime }(n_{2}).$

Consider the equivalence on $M\oplus N$ defined by%
\begin{equation*}
(m,n)\omega (m^{\prime },n^{\prime })\text{ if }g^{\ast }(m)+f^{\ast
}(n)=g^{\ast }(m^{\prime })+f^{\ast }(n^{\prime }).
\end{equation*}%
Clearly, $\omega $ is a congruence. Let%
\begin{equation*}
\pi _{M}^{\omega }:M\longrightarrow (M\oplus N)/\omega ,\pi _{N}^{\omega
}:N\longrightarrow (M\oplus N)/\omega
\end{equation*}%
be the canonical $S$-linear maps, and define 
\begin{eqnarray*}
\pi _{\omega } &:&M\oplus N\rightarrow (M\oplus N)/\omega ,\text{ }%
(m,n)\mapsto \lbrack (m,n)]_{\omega }; \\
h &:&(M\oplus N)/\omega \longrightarrow Q,\text{ }[(m,n)]\mapsto g^{\ast
}(m)+f^{\ast }(n).
\end{eqnarray*}%
Notice that $h$ is well defined by the definition of $\omega $. Then $(\pi
_{M}^{\omega },\pi _{N}^{\omega },(M\oplus N)/\omega )\in \mathcal{P}$.
Since $(g^{\prime },f^{\prime },P)$ is, by assumption, a largest element in $%
\mathcal{P}$, there exists $\alpha :P\rightarrow (M\oplus N)/\omega $ such
that $\alpha \circ \pi _{(g^{\prime },f^{\prime })}=\pi _{\omega }$. It
follows that%
\begin{equation*}
\begin{array}{ccccc}
\varphi (g^{\prime }(m_{1})+f^{\prime }(n_{1})) & = & g^{\ast
}(m_{1})+f^{\ast }(n_{1}) & = & h([(m_{1},n_{1})]_{\omega }) \\ 
& = & (h\circ \pi _{\omega })(m_{1},n_{1}) & = & (\alpha \circ \pi
_{(g^{\prime },f^{\prime })})(m_{1},n_{1}) \\ 
& = & \alpha (g^{\prime }(m_{1})+f^{\prime }(n_{1})) & = & \alpha (g^{\prime
}(m_{2})+f^{\prime }(n_{2})) \\ 
& = & (\alpha \circ \pi _{(g^{\prime },f^{\prime })})(m_{2},n_{2}) & = & 
(h\circ \pi _{\omega })(m_{2},n_{2}) \\ 
& = & h([(m_{2},n_{2})]_{\omega }) & = & g^{\ast }(m_{2})+f^{\ast }(n_{2})
\\ 
& = & \varphi (g^{\prime }(m_{2})+f^{\prime }(n_{2})). &  & 
\end{array}%
\end{equation*}%
Hence $\varphi $ is well defined.$\blacksquare $
\end{itemize}
\end{Beweis}
\end{thm}

\begin{cor}
\label{P-con}Let $f:L\rightarrow M$ and $g:L\rightarrow N$ be morphisms of
left $S$-semimodules. There exists a congruence relation $\rho $ on $M\oplus
N$ such that 
\begin{equation*}
(g^{\prime },f^{\prime };(M\oplus N)/\rho ),\text{ }g^{\prime
}(m):=[(m,0)]_{\rho },\text{ }f^{\prime }(n):=[(0,n)]_{\rho }
\end{equation*}%
is a pushout of $(f,g)$.

\begin{Beweis}
Let $(g^{\ast },f^{\ast },P)$ be a largest element in the poset $(\mathcal{P}%
,\leq )$ in the proof of Theorem \ref{existpushout}. Then $(g^{\ast
},f^{\ast };P)$ is a pushout and there is an surjective map 
\begin{equation*}
\pi :M\oplus N\longrightarrow P,\text{ }(m,n)\mapsto g^{\ast }(m)+f^{\ast
}(n).
\end{equation*}%
Consider the congruence relation $\rho :=\equiv _{\pi }$ and define 
\begin{eqnarray*}
g^{\prime } &:&M\rightarrow (M\oplus N)/\rho ,\text{ }m\mapsto \lbrack
(m,0)]_{\rho } \\
f^{\prime } &:&N\rightarrow (M\oplus N)/\rho ,\text{ }n\mapsto \lbrack
(0,n)]_{\rho }.
\end{eqnarray*}%
For every $l\in L,$ we have%
\begin{equation*}
(g^{\prime }\circ f)(l)=[(f(l),0)]_{\rho }=[(0,g(l))]_{\rho }=(f^{\prime
}\circ g)(l).
\end{equation*}%
The middle equality follows since%
\begin{equation*}
\begin{tabular}{lllll}
$\pi ((f(l),0))$ & $=$ & $(g^{\ast }\circ f)(l)+f^{\ast }(0)$ & $=$ & $%
(f^{\ast }\circ g)(l)+0$ \\ 
& $=$ & $g^{\ast }(0)+(f^{\ast }\circ g)(l)$ & $=$ & $\pi ((0,g(l)))$%
\end{tabular}%
\end{equation*}%
With the canonical map $\pi _{\rho }:M\oplus N\rightarrow (M\oplus N)/\rho ,$
we have $(g^{\prime },f^{\prime },(M\oplus N)/\rho )\in \mathcal{P}$.
Moreover $P\leq (M\oplus N)/\rho $ noticing that%
\begin{equation*}
\alpha :(M\oplus N)/\rho \longrightarrow P,\text{ }[m,n]_{\rho }\mapsto
g^{\ast }(m)+f^{\ast }(n)
\end{equation*}%
is $S$-linear such that $\alpha \pi _{\rho }=\pi _{(g^{\ast },f^{\ast })}.$
Since $(g^{\ast },f^{\ast },P)$ is a largest element in $\mathcal{P},$ $%
(g^{\prime },f^{\prime },(M\oplus N)/\rho )$ is also a largest element in $%
\mathcal{P}.$ Thus $(g^{\prime },f^{\prime };(M\oplus N)/\rho )$ is a
pushout of $(f,g).\blacksquare $
\end{Beweis}
\end{cor}

\begin{lem}
\label{transfers}Let $(g^{\prime },f^{\prime };P)$ be a pushout of the
morphisms of left $S$-semimodules $f:L\rightarrow M$ and $g:L\rightarrow N.$%
%
%
%
%
%
%
%
%
%
%
%
%
%
%
%
%
%
%
%
%
%
%
%
%
%
%
%
%
%
%
%
%
%
%
%
%
%
%
%
%
%
%
%
%
%
%
%
%
%
%
%
%
%
%
%
%
%
%
%
%
%
%
%
%
%
%
%
%
%
%
%
%
%
%
%
%
%

\begin{enumerate}
\item[(1)] If $f$ is surjective, then $f^{\prime }$ is surjective.

\item[(2)] If $f$ is $i$-normal (i.e. $f(L)\subseteq M$ is subtractive),
then $f^{\prime }$ is $i$-normal (i.e. $f^{\prime }(N)\subseteq P$) is
subtractive.

\item[(3)] If $f$ is a normal epimorphism, then $f^{\prime }$ is a normal
epimorphism.

\item[(4)] If $f$ is injective and $g$ is a normal epimorphism, then $%
f^{\prime }$ is injective.
\end{enumerate}

\begin{equation*}
\xymatrix{L \ar[r]^{f} \ar[d]_{g} & M \ar[d]^{g'} \\ N \ar[r]^{f'} & P }
\end{equation*}

\begin{Beweis}
Let $(g^{\prime },f^{\prime };P)$ be a pushout of $(f,g).$

\begin{enumerate}
\item Let $p\in P.$ Since $\pi _{(g^{\prime },f^{\prime })}$ is surjective,
there exists $(m,n)\in M\oplus N$ such that $p=\pi _{(g^{\prime },f^{\prime
})}(m,n)=g^{\prime }(m)+f^{\prime }(n).$ Since $f$ surjective, there exists $%
l\in L$ such that $f(l)=m$. Consider $g(l)+n\in N.$ It follows that%
\begin{equation*}
f^{\prime }(g(l)+n)=(f^{\prime }\circ g)(l)+f^{\prime }(n)=(g^{\prime }\circ
f)(l)+f^{\prime }(n)=g^{\prime }(m)+f^{\prime }(n)=p.
\end{equation*}

\item Let $p\in P$ be such that $p+f^{\prime }(n_{1})=f^{\prime }(n_{2})$
for some $n_{1},n_{2}\in N.$ Pick $(m,n)\in M\oplus N$ such that $p=\pi
_{(g^{\prime },f^{\prime })}(m,n)=g^{\prime }(m)+f^{\prime }(n).$ Thus $%
g^{\prime }(m)+f^{\prime }(n+n_{1})=f^{\prime }(n_{2})$.%
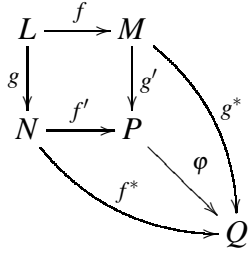
\begin{equation*}
\xymatrix{L \ar[r]^{f} \ar[d]_{g} & M \ar[d]^{g'} \ar@/^1.0pc/[rdd]^{g^*}\\
N \ar[r]^{f'} \ar@/^-1.0pc/[rrd]^{f^*} & P \ar[rd]^{\varphi}\\ & & Q }
\end{equation*}%
Let $\varphi $ be the map from $P$ to the $C$-pushout $Q$ such that $\varphi
\circ g^{\prime }=g^{\ast }$ and $\varphi \circ f^{\prime }=f^{\ast }.$ Then%
\begin{equation*}
\lbrack (m,n+n_{1})]_{\sim }=\varphi (g^{\prime }(m)+f^{\prime
}(n+n_{1}))=\varphi (f^{\prime }(n_{2}))=[(0,n_{2})]_{\sim }.
\end{equation*}%
By the definition of the congruence relation $\sim $ (\ref{CP-cong}), there
exist $l_{1},l_{2}\in L$ such that $m+f(l_{1})=f(l_{2})$ and $%
n+n_{1}+g(l_{2})=n_{2}+g(l_{1})$. Since $f(L)\subseteq M$ is subtractive, $%
m=f(l)$ for some $l\in L$. Then we have%
\begin{equation*}
p=g^{\prime }(m)+f^{\prime }(n)=(g^{\prime }\circ f)(l)+f^{\prime
}(n)=(f^{\prime }\circ g)(l)+f^{\prime }(n)=f^{\prime }(g(l)+n).
\end{equation*}%
It follows that $f^{\prime }(N)\subseteq P$ is subtractive.

\item Without loss of generality, let the pushout be $P=(g^{\prime
},f^{\prime };(M\oplus N)/\rho )$ for some congruence relation $\rho $ on $%
M\oplus N$ and $g^{\prime },f^{\prime }$ are the canonical maps (see
Corollary \ref{P-con}). Since $f$ is surjective, it follows by (1) that $%
f^{\prime }$ is surjective as well.

\textbf{Step I: }Consider the canonical $S$-linear map 
\begin{equation*}
f^{\ast }:N\rightarrow N/Ker(f^{\prime }).
\end{equation*}%
Let $m\in M$ and pick $l\in L$ such that $m=f(l).$ Define 
\begin{equation*}
g^{\ast }:M\rightarrow N/Ker(f^{\prime }),\text{ }m\mapsto (f^{\ast }\circ
g)(l).
\end{equation*}%
\textbf{Claim:} $g^{\ast }$ is well-defined.

Suppose that $f(l)=m=f(l^{\prime })$ for some $l,l^{\prime }\in L.$ Since $f$
is $k$-normal, there exist $l_{1},$ $l_{2}\in Ker(f)$ such that $%
l+l_{1}=l^{\prime }+l_{2}.$ It follows that $g(l)+g(l_{1})=g(l+l_{1})=g(l^{%
\prime }+l_{2})=g(l^{\prime })+g(l_{2})$ with $(f^{\prime }\circ
g)(l_{1})=(g^{\prime }\circ f)(l_{1})=0=(g^{\prime }\circ
f)(l_{2})=(f^{\prime }\circ g)(l_{2}).$ Thus $(f^{\ast }\circ
g)(l)=[g(l)]_{Ker(f^{\prime })}=[g(l^{\prime })]_{Ker(f^{\prime })}=(f^{\ast
}\circ g)(l^{\prime }).$ Clearly, $g^{\ast }$ is $S$-linear and satisfies $%
g^{\ast }\circ f=f^{\ast }\circ g.$

\textbf{Step II: }Define 
\begin{equation*}
\psi :N/Ker(f^{\prime })\rightarrow P,\text{ }[n]_{Ker(f^{\prime })}\mapsto
\lbrack (0,n)]_{\rho }.
\end{equation*}%
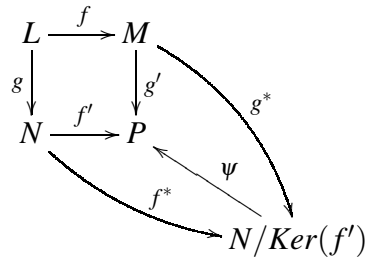
\begin{equation*}
\xymatrix{L \ar[r]^{f} \ar[d]_{g} & M \ar[d]^{g'} \ar@/^1.0pc/[rdd]^{g^*}\\
N \ar[r]^{f'} \ar@/^-1.0pc/[rrd]^{f^*} & P \\ & & N/Ker(f') \ar[lu]_{\psi} }
\end{equation*}%
\textbf{Claim:} $\psi $ is well defined.

Suppose that $[n]_{Ker(f^{\prime })}=[n^{\prime }]_{Ker(f^{\prime })}$ for
some $n,n^{\prime }\in N.$ It follows that $n+n_{1}=n+n_{2}$ for some $%
n_{1},n_{2}\in Ker(f^{\prime })$. Thus%
\begin{equation*}
\begin{array}{ccccc}
\lbrack (0,n)]_{\rho } & = & [(0,n)]_{\rho }+[(0,0)]_{\rho } & = & 
[(0,n)]_{\rho }+f^{\prime }(n_{1}) \\ 
& = & [(0,n)]_{\rho }+[(0,n_{1})]_{\rho } & = & [(0,n+n_{1})]_{\rho } \\ 
& = & [(0,n^{\prime }+n_{2})]_{\rho } & = & [(0,n^{\prime })]_{\rho }.%
\end{array}%
\end{equation*}

For $m\in M,$ pick some $l\in L$ with $f(l)=m.$ Then we have%
\begin{equation*}
\begin{array}{ccccc}
(\psi \circ g^{\ast })(m) & = & (\psi \circ f^{\ast }\circ g)(l) & = & \psi
([g(l)]_{Ker(f^{\prime })}) \\ 
& = & [(0,g(l))]_{\rho } & = & (f^{\prime }\circ g)(l) \\ 
& = & (g^{\prime }\circ f)(l) & = & g^{\prime }(m),%
\end{array}%
\end{equation*}%
whence $\psi \circ g^{\ast }=g^{\prime }.$

On the other hand, for every $n\in N$ we have $(\psi \circ f^{\ast
})(n)=\psi ([n]_{Ker(f^{\prime })})=[(0,n)]_{\rho }=f^{\prime }(n),$ whence $%
(\psi \circ f^{\ast })=f^{\prime }$.

\textbf{Step III: }Since $P$ is a pushout, there exists an $S$-linear map $%
\varphi :P\rightarrow N/Ker(f^{\prime })$ such that $\varphi \circ g^{\prime
}=g^{\ast }$ and $\varphi \circ f^{\prime }=f^{\ast }.$ For each $(m,n)\in
M\oplus N$ we have%
\begin{equation*}
\begin{tabular}{lll}
$(\psi \circ \varphi )([(m,n)]_{\rho })$ & $=$ & $\psi (\varphi
([(m,0)]_{\rho })+\varphi ([(0,n)]_{\rho })$ \\ 
& $=$ & $\psi ((\varphi \circ g^{\prime })(m)+(\varphi \circ f^{\prime
})(n)) $ \\ 
& $=$ & $(\psi \circ g^{\ast })(m)+(\psi \circ f^{\ast })(n))$ \\ 
& $=$ & $g^{\prime }(m)+f^{\prime }(n)$ \\ 
& $=$ & $[(m,n)]_{\rho }.$%
\end{tabular}%
\end{equation*}%
On the other hand, we have for every $n\in N:$%
\begin{equation*}
(\varphi \circ \psi )([n]_{Ker(f^{\prime })})=\varphi \lbrack (0,n)]_{\rho
}=(\varphi \circ f^{\prime })(n)=f^{\ast }(n)=[n]_{Ker(f^{\prime })}.
\end{equation*}%
Hence $P\simeq N/Ker(f^{\prime }).$ This implies that $f^{\prime }$ is $k$%
-normal (as $f^{\ast }$ is obviously $k$-normal).

\item Without loss of generality, let the pushout be $P=(g^{\prime
},f^{\prime };(M\oplus N)/\rho )$ for some congruence relation $\rho $ on $%
M\oplus N$ and $g^{\prime },f^{\prime }$ are the canonical maps (see
Corollary \ref{P-con}). Let $K:=f(Ker(g))$ and consider the canonical
projection $\tilde{g}:M\rightarrow M/K.$ By assumption, $g$ is surjective
and so there exists for every $n\in N$ some $l_{n}\in L$ such that $%
n=g(l_{n}).$

\textbf{Step I:} Define 
\begin{equation*}
\tilde{f}:N\rightarrow M/K,\text{ }n\mapsto \lbrack f(l_{n})]_{K}.
\end{equation*}%
\textbf{Claim:} $\tilde{f}$ is well defined.

Suppose that $g(l_{n})=n=g(l_{n}^{\prime }).$ Since $g$ is $k$-normal, there
exist $l_{1},l_{2}\in Ker(g)$ such that $l_{n}+l_{1}=l_{n}^{\prime }+l_{2},$
whence $f(l_{n})+f(l_{1})=f(l_{n}^{\prime })+f(l_{2}),$ i.e. $%
[f(l_{n})]_{K}=[f(l_{n}^{\prime })]_{K}$ (recall that we chose $K:=f(Ker(g))$%
).%
\begin{equation*}
\xymatrix{L \ar[r]^{f} \ar[d]_{g} & M \ar[d]^{g'}
\ar@/^1.0pc/[rdd]^{\tilde{g}}\\ N \ar[r]^{f'} \ar@/^-1.0pc/[rrd]^{\tilde{f}}
& P \ar[rd]^{\varphi}\\ & & M/K }
\end{equation*}%
Notice that for every $l\in L,$ we have: $(\tilde{f}\circ g)(l)=[f(l)]_{K}=(%
\tilde{g}\circ f)(l)$. Since $P$ is a pushout, there exists an $S$-linear
map $\varphi :P\rightarrow M/K$ such that $(\varphi \circ g^{\prime })=%
\tilde{g}$ and $(\varphi \circ f^{\prime })=\tilde{f}$.

\textbf{Step II:} Define 
\begin{equation*}
\psi :M/K\rightarrow P,\text{ }[m]_{K}\mapsto \lbrack (m,0)]_{\rho }.
\end{equation*}%
We claim that $\psi $ is \emph{well defined}. Suppose that $%
[m]_{K}=[m^{\prime }]_{K}$ for some $m,m^{\prime }\in M.$ Then there exist $%
l_{1},l_{2}\in Ker(g)$ such that $m+f(l_{1})=m^{\prime }+f(l_{2}).$ It
follows that%
\begin{equation*}
\begin{array}{ccccc}
\lbrack (m,0)]_{\rho } & = & g^{\prime }(m) & = & g^{\prime }(m)+0 \\ 
& = & g^{\prime }(m)+(f^{\prime }\circ g)(l_{1}) & = & g^{\prime
}(m)+(g^{\prime }\circ f)(l_{1}) \\ 
& = & g^{\prime }(m+f(l_{1})) & = & g^{\prime }(m^{\prime }+f(l_{2})) \\ 
& = & [m^{\prime },0]_{\rho }. &  & 
\end{array}%
\end{equation*}%
\textbf{Step III:} Notice that for every $n=f(l_{n})\in N$ we have: 
\begin{equation*}
\begin{array}{ccccc}
(\psi \circ \tilde{f})(n) & = & \psi \lbrack f(l_{n})]_{K} & = & 
[(f(l_{n}),0)]_{\rho } \\ 
& = & (g^{\prime }\circ f)(l_{n}) & = & (f^{\prime }\circ g)(l_{n}) \\ 
& = & f^{\prime }(n), &  & 
\end{array}%
\end{equation*}
and%
\begin{equation*}
\begin{array}{ccccc}
(\psi \circ \tilde{g})(m) & = & \psi ([m]_{K}) & = & [(m,0)]_{\rho } \\ 
& = & g^{\prime }(m), &  & 
\end{array}%
\end{equation*}%
thus $\psi \circ \tilde{f}=f^{\prime }$ and $\psi \circ \tilde{g}=g^{\prime
}.$ Moreover, 
\begin{equation*}
\begin{array}{ccccc}
(\varphi \circ \psi )([m]_{K}) & = & \varphi \lbrack (m,0)]_{\rho } & = & 
(\varphi \circ g^{\prime })(m) \\ 
& = & \tilde{g}(m) & = & [m]_{K},\text{ and}%
\end{array}%
\end{equation*}%
\begin{equation*}
\begin{array}{ccccc}
(\psi \circ \varphi )([(m,0)]_{\rho }) & = & (\psi \circ \varphi \circ
g^{\prime })(m) & = & (\psi \circ \tilde{g})(m) \\ 
& = & \psi ([m]_{K}) & = & [(m,0)]_{\rho },%
\end{array}%
\end{equation*}%
i.e. $\psi ,\varphi $ are $S$-linear isomorphisms and $\psi ^{-1}=\varphi .$
Moreover, $M/K$ is a pushout.

\textbf{Step IV:} Let $n,n^{\prime }\in N$ be such that $\tilde{f}(n)=\tilde{%
f}(n^{\prime }),$ i.e. $[f(l_{n})]_{K}=[f(l_{n^{\prime }})]_{K}.$ Then there
exist $l_{1},l_{2}\in Ker(g)$ such that $f(l_{n}+l_{1})=$ $%
f(l_{n})+f(l_{1})=f(l_{n^{\prime }})+f(l_{2})=f(l_{n^{\prime }}+l_{2}),$
whence $l_{n}+l_{1}=l_{n^{\prime }}+l_{2}$ as $f$ is injective. It follows
that $n=g(l_{n})=g(l_{n})+g(l_{1})=g(l_{n}+l_{1})=g(l_{n^{\prime
}}+l_{2})=g(l_{n^{\prime }})+g(l_{2})=g(l_{n^{\prime }})=n^{\prime }$. Thus $%
\tilde{f}$ is injective. Since $f^{\prime }=\psi \circ \tilde{f}$ and $\psi ,%
\tilde{f}$ are injective, we conclude that $f^{\prime }$ is injective as
well.$\blacksquare $
\end{enumerate}
\end{Beweis}
\end{lem}

\section{Projective Semimodules}

\label{PIF}

As before, $(S,+,0,\cdot ,1)$ is a semiring and, unless otherwise explicitly
mentioned, an $S$-module is a \textbf{left} $S$-semimodule. Exact sequences
here are in the sense of Abuhlail \cite{Abu2014} (Definition \ref{Abu-exs}).

\markright{\scriptsize\tt Chapter \ref{sec-proj}: $e$-projective
Semimodules} 

There are several notions of projectivity for a semimodule over a semiring,
which coincide if it were a module over a ring. In this Chapter, we consider
some of them and clarify the relationships between them, and then
investigate the so called $e$\emph{-projective semimodules} which turn to
coincide with the so called \emph{normally projective semimodules} (both
notions introduced by Abuhlail \cite[1.25, 1.24]{Abu2014-CA} and called 
\emph{uniformly projective semimodules}). The terminology ``$e$\emph{%
-projective}" appeared first in \cite{AIKN2018}).

\begin{definition}
(\cite{AIKN2018})\ A left $S$-semimodules $P$ is

$M$-$e$-\textbf{projective}%
\index{Semimodule!$e$-projective} (where $M$ is a left $S$-semimodule) if
the covariant functor 
\begin{equation*}
Hom_{S}(P,-):%
\text{ }_{S}\mathbf{SM}\longrightarrow \text{ }_{\mathbb{Z}^{+}}\mathbf{SM}
\end{equation*}%
transfers every short exact sequence of left $S$-semimodules%
\begin{equation}
0\longrightarrow L\overset{f}{\longrightarrow }M\overset{g}{\longrightarrow }%
N\longrightarrow 0  \label{lmn-k}
\end{equation}%
into a short exact sequence of commutative monoids%
\begin{equation}
0\longrightarrow Hom_{S}(P,L)\overset{(P,f)}{\longrightarrow }Hom_{S}(P,M)%
\overset{(P,g)}{\longrightarrow }Hom_{S}(P,N)\longrightarrow 0.
\label{hom(p,lmn)}
\end{equation}%
We say that $P$ is $e$\textbf{-projective}%
\index{e-projective} if $P$ is $M$-$e$-projective for every left $S$%
-semimodule $M$.
\end{definition}

\begin{punto}
Let $P$ be a left $S$-semimodule.

For a left $S$-semimodule $M,$ we say that $P$ is

$M$\textbf{-projective}%
\index{Semimodule!projective} \cite[page 195]{Gol1999} if for every \emph{%
surjective} $S$-linear map $f:M\rightarrow N$ and an $S$-linear map $%
g:P\rightarrow N,$ there exists an $S$-linear map $h:P\rightarrow M$ such
that $f\circ h=g;$%
\begin{equation*}
\xymatrix{ M \ar[r]^{f} & N \ar[r] & 0 \\ & P \ar[u]_{g} \ar@{.>}[lu]^{h} }
\end{equation*}

$M$-$k$\textbf{-projective}%
\index{Semimodule!$k$-projective} \cite[Definition 6]{Alt1996} if for every 
\emph{normal epimorphism} $f:M\rightarrow N$ and any $S$-linear map $%
g:P\rightarrow N,$ there exists an $S$-linear map $h:P\rightarrow M$ such
that $f\circ h=g;$%
\begin{equation*}
\xymatrix{ M \ar[rrr]^{f (normal)} & & & N \ar[r] & 0\\ & & & P \ar[u]_{g}
\ar@{.>}[lllu]^{h} }
\end{equation*}

\textbf{normally }$M$-\textbf{projective }%
\index{Semimodule!normally projective}\cite[1.25]{Abu2014-CA} if for every 
\emph{normal epimorphism} $f:M\rightarrow N$ and any $S$-linear map $%
g:P\rightarrow N,$ there exists an $S$-linear map $h:P\rightarrow M$ such
that $f\circ h=g$%
\begin{equation*}
\xymatrix{P \ar@<-.5ex>[r]_{h_1} \ar@<.5ex>[r]^{h_2} & M \ar[rrr]^{f
(normal)} & & & N \ar[r] & 0\\ \\ & & & & P \ar[uu]_{g} \ar@{.>}[llluu]^{h}
\ar@/^-0.5pc/[llluu]_{h'}}
\end{equation*}%
and whenever an $S$-linear map $h^{\prime }:P\rightarrow M$ satisfies $%
f\circ h^{\prime }=g$, there exist $S$-linear maps $h_{1},h_{2}:P\rightarrow
M$ such that $f\circ h_{1}=0=f\circ h_{2}$ and $h+h_{1}=h^{\prime }+h_{2}$.

We say that $P$ is \textbf{projective} (resp., $k$\textbf{-projective}, 
\textbf{normally projective}) if $P$ is $M$-projective (resp., $M$-$k$%
-projective, normally $M$-projective) for every left $S$-semimodule $M.$
\end{punto}

\begin{prop}
\label{retract}\emph{(}cf., \emph{\cite[Theorem 1.9]{Tak1983}, \cite[%
Proposition 17.16]{Gol1999}}\emph{)} A left $S$-semimodule $_{S}P$ is
projective if and only if $P$ is a retract of a free left $S$-semimodule.
\end{prop}

\begin{rems}
\begin{enumerate}
\item It is obvious that projective and $e$-projective semimodules are $k$%
-projective.

\item Despite being a retract of a free semimodule, a projective semimodule
is not necessarily a direct summand of a free semimodule (\cite[Example 2.3]%
{Alt2002}).
\end{enumerate}
\end{rems}

\begin{prop}
\label{e=n}Let $P$ be a left $S$-semimodule.

\begin{enumerate}
\item Let $M$ be a left $S$-semimodule. Then $_{S}P$ is $M$-$e$-projective
if and only if $_{S}P$ is normally $M$-projective.

\item $_{S}P$ is $e$-projective if and only if $_{S}P$ is normally
projective.
\end{enumerate}
\end{prop}

\begin{Beweis}
We need to prove (1) only.

($\Longrightarrow $) Assume that $_{S}P$ is $M$-$e$-projective. Let $%
f:M\rightarrow N$ be a normal epimorphism and $g:P\rightarrow N$ an $S$%
-linear map. By Lemma \ref{exact}, the sequence%
\begin{equation*}
0\longrightarrow Ker(f)\overset{\iota }{\longrightarrow }M\overset{f}{%
\longrightarrow }N\longrightarrow 0
\end{equation*}%
is a short exact sequence, where $\iota $ is the canonical embedding. By
assumption, the following sequence of commutative monoids%
\begin{equation*}
0\longrightarrow Hom_{S}(P,Ker(f))\overset{(P,\iota )}{\longrightarrow }%
Hom_{S}(P,M)\overset{(P,g)}{\longrightarrow }Hom_{S}(P,N)\longrightarrow 0
\end{equation*}%
is exact. In particular, $(P,g)$ is surjective and $k$-normal, whence $P$ is
normally $M$-projective.

($\Longleftarrow $) let $0\longrightarrow L\overset{f}{\longrightarrow }M%
\overset{g}{\longrightarrow }N\longrightarrow 0$ be a short exact sequence
of left $S$-semimodules and consider the induces sequences of commutative
monoids%
\begin{equation*}
0\longrightarrow Hom_{S}(P,L)\overset{(P,f)}{\longrightarrow }Hom_{S}(P,M)%
\overset{(P,g)}{\longrightarrow }Hom_{S}(P,N)\longrightarrow 0.
\end{equation*}%
By Proposition \ref{lr-exact}, $(P,f)$ is a normal monomorphism and $%
\func{Im}((P,f))=Ker((P,g)).$ By assumption, $(P,g)$ is a normal
epimorphism, whence the induced sequence of commutative monoids is exact.$%
\blacksquare $
\end{Beweis}

Following an observation by H. Al-Thani made in \cite[theorem 4]{Alt1995},
we provide a \emph{detailed proof} that every projective $S$-semimodule is $%
e $-projective.

\begin{prop}
\label{proj->uproj}Every projective left $S$-semimodule is $e$-projective.
\end{prop}

\begin{Beweis}
Let $_{S}P$ be projective. Assume that $M\overset{g}{\longrightarrow }%
N\longrightarrow 0$ is a normal epimorphism of left $S$-semimodules, and $%
\alpha \in \mathrm{Hom}_{S}(P,N).$ Since $_{S}P$ is $M$-projective,%
\begin{equation*}
\mathrm{Hom}_{S}(P,M)\overset{(P,g)}{\longrightarrow }\mathrm{Hom}%
_{S}(P,N)\longrightarrow 0
\end{equation*}%
is surjective, \emph{i.e.} there exists $\beta \in \mathrm{Hom}_{S}(P,M)$
such that $g\circ \beta =\alpha .$

By Proposition \ref{e=n}, it is enough to prove that $(P,g)$ is $k$-normal.

Suppose that $(P,g)(\beta )=(P,g)(\beta ^{\prime })$ for some $\beta ,\beta
^{\prime }\in \mathrm{Hom}_{S}(P,M),$ i.e. $g\circ \beta =g\circ \beta
^{\prime }.$ Since $_{S}P$ is projective, $P$ is a retract of a free left $S$%
-semimodule, \emph{i.e.} there exists an index set $\Lambda $ and a
surjective $S$-linear map $\theta :S^{(\Lambda )}\longrightarrow P$ as well
as an injective $S$-linear map $\psi :P\longrightarrow S^{(\Lambda )}$ such
that $\theta \circ \psi =id_{P}.$ Notice that $g\circ \beta \circ \theta
=g\circ \beta ^{\prime }\circ \theta .$ For every $\lambda \in \Lambda ,$
and since $g$ is $k$-normal, there exist $m_{\lambda },$ $m_{\lambda
}^{\prime }\in \mathrm{Ker}(g)$ such that $(\beta \circ \theta )(\lambda
)+m_{\lambda }=(\beta ^{\prime }\circ \theta )(\lambda )+m_{\lambda
}^{\prime }.$ Let $\gamma ,\gamma ^{\prime }\in \mathrm{Hom}_{S}(S^{(\Lambda
)},M)$ be the \emph{unique} $S$-linear maps with $\gamma (\lambda
)=m_{\lambda }$ and $\gamma ^{\prime }(\lambda )=m_{\lambda }^{\prime }$ for
each $\lambda \in \Lambda $ (they exist and are unique since $\Lambda $ is a
basis for $S^{(\Lambda )}$). It follows that 
\begin{equation*}
g\circ (\gamma \circ \psi )=(g\circ \gamma )\circ \psi =0=(g\circ \gamma
^{\prime })\circ \psi =g\circ (\gamma ^{\prime }\circ \psi ),
\end{equation*}%
i.e. $\gamma \circ \psi ,$ $\gamma ^{\prime }\circ \psi \in \mathrm{Ker}%
((P,g)).$ Moreover, for any $\lambda \in \Lambda $ we have%
\begin{equation*}
(\beta \circ \theta +\gamma )(\lambda )=(\beta \circ \theta )(\lambda
)+m_{\lambda }=(\beta ^{\prime }\circ \theta )(\lambda)+m_{\lambda }^{\prime
}=(\beta ^{\prime }\circ \theta +\gamma ^{\prime })(\lambda ),
\end{equation*}%
whence $\beta \circ \theta +\gamma =\beta ^{\prime }\circ \theta +\gamma
^{\prime }.$ It follows that%
\begin{equation*}
\begin{array}{ccccc}
\beta +\gamma \circ \psi & = & \beta \circ id_{P}+\gamma \circ \psi & = & 
\beta \circ (\theta \circ \psi )+\gamma \circ \psi \\ 
& = & (\beta \circ \theta +\gamma )\circ \psi & = & (\beta ^{\prime }\circ
\theta +\gamma ^{\prime })\circ \psi \\ 
& = & \beta ^{\prime }\circ (\theta \circ \psi )+\gamma ^{\prime }\circ \psi
& = & \beta ^{\prime }\circ id_{P}+\gamma ^{\prime }\circ \psi \\ 
& = & \beta ^{\prime }+\gamma ^{\prime }\circ \psi .\blacksquare &  & 
\end{array}%
\end{equation*}
\end{Beweis}

The following example shows that the class of $S$-$e$-projective left $S$%
-semimodules is strictly larger than that of $S$-projective left $S$%
-semimodules.

\begin{ex}
\label{e-proj-not-proj}Consider the semiring $S:=\mathbb{Q}^{+}$ of
non-negative rational numbers, with the usual addition and multiplication.
Consider the Boolean algebra $\mathbb{B}=\{0,1\}$ as an $S$-semimodule with $%
s\cdot 1=1\Leftrightarrow s\in S\backslash \{0\}$. Then $_{S}\mathbb{B}$ is $%
S$-$e$-projective but \emph{not} $S$-projective.
\end{ex}

\begin{Beweis}
Consider the $S$-linear map%
\begin{equation*}
f:S\rightarrow \mathbb{B},\text{ }s\mapsto \left\{ 
\begin{array}{ccc}
1 & , & s\neq 0 \\ 
&  &  \\ 
0, &  & s=0%
\end{array}%
\right\vert
\end{equation*}
Notice that $f$ is not $k$-normal: $Ker(f)=\{0\}$, $f(1)=1=f(2)$, and $%
1+0\neq 2+0$.

Since there is no surjective $S$-linear map from $\mathbb{B}$ to $S$, there
is no isomorphism from $\mathbb{B}$ to $S$. Since $S$ is an ideal-simple $S$%
-semimodule, $Hom_{S}(\mathbb{B},S)=\{0\}$ by Lemma \ref{s-char}. Since the
following diagram%
\begin{equation*}
\xymatrix{S \ar[r]^{f} & \mathbb{B} \\ & \mathbb{B} \ar[u]_{id_{\mathbb{B}}}
\ar[ul]^0 }
\end{equation*}%
cannot be completed commutatively, $B$ is not $S$-projective.

Let $N$ be an $S$-semimodule and $f:S\rightarrow N$ be a normal $S$%
-epimorphism. If $f=0$, then $N=f(S)=0$, which implies that every $S$-linear
map $g:\mathbb{B}\rightarrow N$ is the zero morphism and by choosing $S$%
-linear map $0=h:\mathbb{B}\rightarrow S$ we have $g=f\circ h$.

If $f\neq 0$, then $f(1)\neq 0$. For every $s\in S\backslash \{0\}$, we have 
$0\neq f(1)=f(s^{-1}s)=s^{-1}f(s)$, whence $f(s)\neq 0$. Thus $Ker(f)=\{0\}$%
. If $f(s)=f(t)$, then $s+k_{1}=t+k_{2}$ for some $k_{1},k_{2}\in
Ker(f)=\{0\}$, thus $s=t$. Hence, $f$ is an $S$-isomorphism. Since $S$ is
not $S$-isomorphic to $\mathbb{B}$, $N$ is not $S$-isomorphic to $\mathbb{B}$%
. Since $S$ is ideal-simple, $N$ is ideal-simple. Thus $Hom_{S}(\mathbb{B}%
,N)=\{0\}$ and $\mathbb{B}$ is $S$-$e$-projective.$\blacksquare $
\end{Beweis}

\begin{prop}
\label{lmn-prop}Let%
\begin{equation}
L\overset{f}{\longrightarrow }M\overset{g}{\longrightarrow }N  \label{LMN-3}
\end{equation}%
be a sequence of left $S$-semimodules, $P$ a left $S$-semimodule and
consider the sequence%
\begin{equation}
Hom_{S}(P,L)\overset{(P,f)}{\longrightarrow }Hom_{S}(P,M)\overset{(P,g)}{%
\longrightarrow }Hom_{S}(P,N)  \label{P-LMN-3}
\end{equation}%
of commutative monoids.

\begin{enumerate}
\item If (\ref{LMN-3}) is exact with $f$ normal and $P$ is $e$-projective,
then (\ref{P-LMN-3}) is exact and $(P,f)$ is normal.

\item If (\ref{LMN-3}) is exact with $f$ normal and $P$ is $k$-projective,
then (\ref{P-LMN-3}) is proper-exact.

\item If (\ref{LMN-3}) is exact and $P$ is projective, then (\ref{P-LMN-3})
is proper-exact.
\end{enumerate}
\end{prop}

\begin{Beweis}
Consider the exact sequence of left $S$-semimodules%
\begin{equation*}
0\longrightarrow Ker(g)\overset{\iota }{\longrightarrow }M\overset{\pi }{%
\longrightarrow }M/Ker(g)\longrightarrow 0
\end{equation*}%
with canonical $S$-linear maps (see Corollary \ref{M/L}). Assume (\ref{LMN-3}%
) to be exact, so that $f(M)=Ker(g)$ and $M/Ker(g)=M/f(M)\simeq Coker(f).$
By the \emph{Universal Property of Kernels}, there exists a unique $S$%
-linear map $\widetilde{f}:L\longrightarrow Ker(g)$ such that $\iota \circ 
\widetilde{f}=f$. On the other hand, by the \emph{Universal Property of
Cokernels, }there exists a unique $S$-linear map $\widetilde{g}%
:M/Ker(g)\longrightarrow N$ such that $\widetilde{g}\circ \pi =g$. So, we
have a commutative diagram of left $S$-semimodules%
\begin{equation}
\xymatrix{ & & & L \ar[lld]_{\widetilde{f}} \ar[d]^{f} & & & 0 \ar[ld] \\ 0
\ar[r] & Ker(g) \ar[ld] \ar[rr]^{\iota} & & M \ar[rr]^{\pi} \ar[d]^{g} & &
M/Ker(g) \ar[lld]^{\widetilde{g}} \ar[r] & 0 \\ 0 & & & N & & & }
\label{M-ker}
\end{equation}%
Applying the contravariant functor $Hom_{S}(P,-),$ we get the sequence%
\begin{equation}
0\longrightarrow Hom_{S}(P,Ker(g))\overset{(P,\iota )}{\longrightarrow }%
Hom_{S}(P,M)\overset{(P,\pi )}{\longrightarrow }Hom_{S}(P,M/Ker(g))%
\longrightarrow 0  \label{Hom-Ker}
\end{equation}%
and we obtain the commutative diagram%
\begin{equation}
\xymatrix{& & & Hom_S(P,L) \ar[dd]^{(P,f)} \ar[lldd]_{(P,{\widetilde{f}})} &
& & 0 \ar[ldd] \\ & & & & & & & \\ 0 \ar[r] & Hom_{S}(P,Ker(g))
\ar[rr]^{(P,\iota)} \ar[ldd] & & Hom_{S}(P,M) \ar[rr]^{(P,\pi)}
\ar[dd]^{(P,g)} & & Hom_{S}(P,M/Ker(g)) \ar[lldd]^{{(P,{\widetilde{g}})}}
\ar[r] & 0 \\ & & & & & & & \\ 0 & & & Hom_S(P,N) & & & &}  \label{homM-ker}
\end{equation}%
of commutative monoids.

Notice that $\widetilde{g}$ is injective since $g$ is $k$-normal, whence $(P,%
\widetilde{g})$ is injective. On the other hand, $\widetilde{f}$ is
surjective since $f(M)=Ker(g)$. If, moreover, $f=\iota \circ \widetilde{f}$
is $k$-normal then it follow by Lemma \ref{i-normal} (1-a) that $\widetilde{f%
}$ is $k$-normal (whence normal).

\begin{enumerate}
\item Let $_{S}P$ be $e$-projective, so that $(P,\widetilde{f})$ is
surjective. It follows then by Proposition \ref{lr-exact} that Sequence (\ref%
{Hom-Ker}) is (proper-)exact.

\textbf{Step I:}\ We have%
\begin{equation*}
\begin{array}{ccccc}
Ker((P,g)) & = & Ker((P,\widetilde{g})\circ (P,\pi )) &  &  \\ 
& = & Ker((P,\pi )) &  & \text{(}(P,\widetilde{g})\text{ is injective)} \\ 
& = & im((P,\iota )) &  & \text{(Proposition \ref{lr-exact} (2))} \\ 
& = & im((P,\iota )\circ (P,\widetilde{f})) &  & \text{(}(P,\widetilde{f})%
\text{ is surjective)} \\ 
& = & im(P,f). &  & 
\end{array}%
\end{equation*}

\textbf{Step II:}\ Since $(P,\widetilde{g})$ is injective and $(P,\pi )$ is $%
k$-normal, it follows by Lemma \ref{i-normal} (1-a) that $(P,g)=(P,%
\widetilde{g})\circ (P,\pi )$ is $k$-normal. So, (\ref{P-LMN-3}) is exact.
Moreover, $(P,\iota )$ is injective and $(P,\widetilde{f})$ is normal,
whence $(P,f)=(P,\iota )\circ (P,\widetilde{f})$ is normal by Lemma \ref%
{i-normal} (2-c).

\item The proof is Step I of (1) noticing that $(P,\widetilde{f})$ is
surjective since $\widetilde{f}$ is normal and $_{S}P$ is $k$-projective.

\item The proof is Step I of (1) noticing that $(P,\widetilde{f})$ is
surjective since $_{S}P$ is projective without assuming that $\widetilde{f}$
is $k$-normal.$\blacksquare $
\end{enumerate}
\end{Beweis}

\begin{thm}
\label{e-proj-3-char}Let $M$ be a left $S$-semimodule. The following are
equivalent for a left $S$-semimodule $P:$

\begin{enumerate}
\item $_{S}P$ is normally $M$-projective;

\item $_{S}P$ is $M$-$e$-projective;

\item For every \emph{exact} sequence of left $S$-semimodules (\ref{LMN-3}),
the included sequence (\ref{P-LMN-3}) of commutative monoids is \emph{exact}
and $(P,f)$ is normal.
\end{enumerate}
\end{thm}

\begin{Beweis}
$(1)\Longleftrightarrow (2)$ follows by Proposition \ref{e=n}.

$(2)\Rightarrow (3)$ follows by Proposition \ref{lmn-prop} (1).

$(3)\Rightarrow (1)$ This follows directly by applying the assumption to the
exact sequences of the form $M\overset{f}{\longrightarrow }N\longrightarrow
0 $ with $f$ normal.$\blacksquare $
\end{Beweis}

Using Propositions \ref{lr-exact} and \ref{lmn-prop}, we recover the
following characterizations of $k$-projective semimodules \cite[Theorem 8]%
{Alt1996} and \cite[Theorem 3.7]{Alt2002}:

\begin{thm}
\label{k-proj-3-char}Let $M$ be a left $S$-semimodule. The following are
equivalent for a left $S$-semimodule $P:$

\begin{enumerate}
\item $_{S}P$ is $M$-$k$-projective;

\item For every \emph{exact} sequence of left $S$-semimodules (\ref{lmn-k}),
the induced sequence (\ref{hom(p,lmn)}) of commutative monoids is \emph{%
proper-exact}.

\item For every \emph{exact} sequence of left $S$-semimodules (\ref{LMN-3})
in which $f$ is normal, the induced sequence (\ref{P-LMN-3}) of commutative
monoids is \emph{proper-exact}.
\end{enumerate}
\end{thm}

Using Propositions \ref{lr-exact} and \ref{lmn-prop}, we recover the
following characterizations of projective semimodules \cite[Theorem 3.5]%
{Alt2002}:

\begin{thm}
\label{proj-3-char}Let $M$ be a left $S$-semimodule. The following are
equivalent for a left $S$-semimodule $P:$

\begin{enumerate}
\item $_{S}P$ is $M$-projective;

\item For every \emph{proper-exact} sequence of left $S$-semimodules (\ref%
{lmn-k}) in which $f$ is normal, the induced sequence (\ref{hom(p,lmn)}) of
commutative monoids is \emph{proper-exact}.

\item For every \emph{exact} sequence of left $S$-semimodules (\ref{LMN-3}),
the induced sequence (\ref{P-LMN-3}) of commutative monoids is \emph{%
proper-exact}.
\end{enumerate}
\end{thm}

\begin{punto}
We call a short exact sequence of $S$-semimodules%
\begin{equation}
0\rightarrow A\overset{f}{\longrightarrow }B\overset{g}{\longrightarrow }%
C\rightarrow 0  \label{se}
\end{equation}

\emph{left splitting} if there exists $f^{\prime }\in \mathrm{Hom}_{S}(B,A)$
such that $f^{\prime }\circ f=id_{A};$

\emph{right splitting} if there exists $g^{\prime }\in \mathrm{Hom}_{S}(C,B)$
such that $g\circ g^{\prime }=id_{C}.$

We say that (\ref{se}) splits or is splitting if it is left splitting and
right splitting.
\end{punto}

Left splitting of short exact sequences of semimodules is not equivalent to
right splitting.

\begin{ex}
Consider the semiring $B(3,1)=(\{0,1,2\},\oplus ,0,\otimes ,1)$, where%
\begin{equation*}
1\oplus 2=1,\text{ }2\oplus 2=0,\text{ }2\otimes 2=0;
\end{equation*}%
see \cite[Example 1.8]{Gol1999}. Then we have a short exact sequence of
commutative monoids%
\begin{equation}
0\longrightarrow \{0,2\}\overset{\iota }{\longrightarrow }B(3,1)\overset{\pi 
}{\longrightarrow }\mathbb{Z}_{2}\longrightarrow 0,  \label{0,p}
\end{equation}%
where $\iota $ is the canonical embedding and $\pi $ is the canonical
projection. The sequence (\ref{0,p}) is exact since $\{0,2\}$ is subtractive
and $B(3,1)/\{0,2\}\cong $ $_{\mathbb{Z}^{+}}\mathbb{Z}_{2}$ (see Lemma \ref%
{exact}). Consider 
\begin{equation*}
f:B(3,1)\longrightarrow \{0,2\},\text{ }x\mapsto \left\{ 
\begin{array}{ccc}
2, &  & x\neq 0 \\ 
&  &  \\ 
0, &  & x=0%
\end{array}%
\right.
\end{equation*}%
and notice that $f\circ \iota =id_{\{0,2\}},$ i.e. (\ref{0,p}) is left
splitting. On the other hand, $Hom_{\mathbb{Z}^{+}}(\mathbb{Z}%
_{2},B(3,1))=\{0\}$. Consequently, (\ref{0,p}) is not right splitting.
\end{ex}

\begin{prop}
\label{char-k-proj}A left $S$-semimodule $_{S}P$ is $k$-projective if and
only if every short exact sequence of left $S$-semimodules 
\begin{equation*}
0\rightarrow A\overset{f}{\longrightarrow }B\overset{g}{\longrightarrow }%
P\rightarrow 0
\end{equation*}%
is right-splitting.\newline
\end{prop}

\begin{Beweis}
($\Rightarrow $) Let $P$ be $k$-projective and $0\rightarrow L\overset{f}{%
\longrightarrow }M\overset{g}{\longrightarrow }P\rightarrow 0$ be a short
exact sequence. In particular, $g$ is surjective and $k$-normal. Consider, $%
id_{P}:P\longrightarrow P.$ Since $_{S}P$ is $k$-projective, there exists an 
$S$-linear map $g^{\prime }:P\rightarrow M$ such that the following diagram%
\begin{equation*}
\xymatrix{& & P \ar[dd]^{id_P} \ar@{.>}[ddll]_{g'} \\ \\ M \ar[rr]_{g} & & P}
\end{equation*}%
is commutative, i.e. $g\circ g^{\prime }=id_{P}.$

($\Leftarrow $) Let $M\overset{g}{\longrightarrow }N\longrightarrow 0$ be a
normal surjective $S$-linear map and $h:P\rightarrow N$ be a morphism of
left $S$-semimodules. Consider the pullback of $g$ and $h:$%
\begin{equation*}
Q:=\{(p,m)\in P\times M\text{ }|\text{ }h(p)=g(m)\}
\end{equation*}%
and the following commutative diagram%
\begin{equation*}
\xymatrix{Q \ar[rr]^{\pi_P} \ar[dd]_{\pi_M} & & P \ar[dd]^{h} \\ \\ M
\ar[rr]_{g} & & N}
\end{equation*}%
where $\pi _{P}$ and $\pi _{Q}$ are the canonical projections. Since $g$ is
surjective, $h(p)=g(m)$ for some $m\in M,$ i.e. $(p,m)\in Q$ and indeed, $%
p=\pi _{P}(p,m).$ Hence $\pi _{P}$ is surjective. Let $(p,m),(p,m^{\prime
})\in Q$ so that $\pi _{P}(p,m)=\pi _{P}(p,m^{\prime }).$ Then $%
g(m)=h(p)=g(m^{\prime })$ and there exist $u,u^{\prime }\in Ker(g)$ such
that $m+u=m^{\prime }+u^{\prime }$ (since $g$ is $k$-normal). Notice that $%
(0,u),(0,v)\in Ker(\pi _{P})$ and $(p,m)+(0,u)=(p,m+u)=(p,m^{\prime
}+u^{\prime })=(p,m)+(0+u^{\prime }),$ i.e. $\pi _{P}$ is $k$-normal. Hence
the sequence 
\begin{equation*}
0\rightarrow Ker(\pi _{P})\hookrightarrow Q\overset{\pi _{P}}{%
\longrightarrow }P\rightarrow 0
\end{equation*}%
is exact, and there exists by our assumption an $S$-linear map $\varphi
:P\rightarrow Q$ such that $\pi _{P}\circ \varphi =id_{P}.$ Notice that for
every $p\in P,$ $\varphi (p)\in Q,$ whence $\varphi (p)=(p,m)$ for some $%
m\in M$ with $h(p)=g(m).$ It follows that 
\begin{equation}
(g\circ (\pi _{M}\circ \varphi ))(p)=g(\pi _{M}(p,m))=g(m)=h(p).
\end{equation}%
\newline
So, $g\circ (\pi _{M}\circ \varphi )=h.$ Consequently, $P$ is $k$-projective.%
$\blacksquare $
\end{Beweis}

\begin{lem}
\label{sumproj}If $M$ is a left $S$-semimodule such that every subtractive
subsemimodule is a direct summand, then every left $S$-semimodule is $M$-$e$%
-projective.

\begin{Beweis}
Let $P$ be a left $S$-semimodule and let 
\begin{equation*}
f:M\longrightarrow N\longrightarrow 0
\end{equation*}%
be a normal epimorphism and $g:P\rightarrow N$ be an $S$-linear map. Notice
that $Ker(f)\leq _{S}M$ is a subtractive subsemimodule, whence $%
M=Ker(f)\oplus L$ for some subsemimodule $L\leq _{S}M.$ The row of this
following diagram is exact by Lemma \ref{exact}%
\begin{equation*}
\xymatrix{0 \ar[r] & Ker(f) \ar[r]^{\iota} & M \ar[r]^{f} & N \ar[r] & 0\\ &
& & P \ar[u]_{g} }
\end{equation*}%
It follows (see also Remark \ref{d-iso}(2)) that we have isomorphisms of
left $S$-semimodules:%
\begin{equation*}
N\simeq M/Ker(f)\simeq L.
\end{equation*}%
Considering the induced isomorphism $N\overset{g^{\prime }}{\simeq }L$ and
setting $h:=\iota _{L}\circ g^{\prime }\circ g:P\longrightarrow M$ where $%
f\circ \iota _{L}=id_{L}$ and $\iota _{L}\circ f|_{L}=id_{N}$, we have
indeed $f\circ h=g.$

Suppose that also $h^{\prime }:P\longrightarrow M$ satisfies $f\circ
h^{\prime }=g.$ Consider the projection $\pi :M\longrightarrow Ker(f).$ Then 
$\varphi :=\iota _{L}\circ g^{\prime }\circ f+\pi =id_{M}:$ Let $m\in M,$
and write $m=k+l$ for some unique $k\in Ker(f)$ and $l\in L,$ and notice that%
\begin{eqnarray*}
\varphi (m) &=&\varphi (k+l) \\
&=&(\iota _{L}\circ g^{\prime }\circ f+\pi )(k+l)+(\iota _{L}\circ g^{\prime
}\circ f)(k+l)+\pi (k+l) \\
&=&l+k \\
&=&m.
\end{eqnarray*}%
\ 

Choose $h_{1}:=\pi \circ h^{\prime }:P\longrightarrow M$ and $%
h_{2}=0:P\longrightarrow M.$ Notice that $f\circ h_{1}=f\circ \pi \circ
h^{\prime }=0=f\circ h_{2}.$ Moreover, we have for each $p\in P:$%
\begin{equation*}
\begin{array}{ccccc}
(h+h_{1})(p) & = & h(p)+h_{1}(p) & = & (\iota _{L}\circ g^{\prime }\circ
g)(p)+(\pi \circ h^{\prime })(p) \\ 
& = & (\iota _{L}\circ g^{\prime }\circ f\circ h^{\prime })(p)+\pi \circ
h^{\prime }(p) & = & ((\iota _{L}\circ g^{\prime }\circ f+\pi )\circ
h^{\prime })(p) \\ 
& = & h^{\prime }(p) & = & (h^{\prime }+0)(p).%
\end{array}%
\end{equation*}%
Consequently, $P$ is $M$-$e$-projective.$\blacksquare $
\end{Beweis}
\end{lem}

The following two results are relative versions of parts of \cite[Corollary
3.3]{AIKN2018}; moreover, we give \emph{detailed homological} proofs as the
ones in \cite{AIKN2018} are compact and categorical.

\begin{lem}
\label{ret-proj}\emph{(cf., \cite[Corollary 3.3]{AIKN2018})}

\begin{enumerate}
\item Let $M$ be a left $S$-semimodule. A retract of an $M$-$e$-projective
semimodule is $M$-$e$-projective.

\item A retract of an $e$-projective left $S$-semimodule is $e$-projective.
\end{enumerate}

\begin{Beweis}
We only need to prove (1).

Let $P$ be a left $S$-semimodule which is $M$-$e$-projective and let $_{S}K$
be a retract of $P$ along with a surjective $S$-linear map $\pi
_{K}:P\rightarrow K$ and an injective $S$-linear map $\iota
_{K}:K\rightarrow P$ such that $\pi _{K}\circ \iota _{K}=id_{K}.$

Let $f:M\rightarrow N$ be a normal epimorphism and $g:K\rightarrow N$ an $S$%
-linear map.

Since $P$ is $e$-projective, there exists an $S$-linear map $h^{\ast
}:P\rightarrow M$ such that $f\circ h^{\ast }=g\circ \pi _{K}$. Consider $%
h:=h^{\ast }\circ \iota _{K}:K\rightarrow M$.%
\begin{equation*}
\xymatrix{ M \ar[r]^{f} & N \ar[r] & 0\\ & K \ar[u]_{g}
\ar@<-.5ex>[d]_{\iota_K} \\ & P \ar@<-.5ex>[u]_{\pi_K}
\ar@/^1.0pc/@{-->}[luu]^{h^{*}}}
\end{equation*}%
Then $f\circ h=f\circ (h^{\ast }\circ \iota _{K})=g\circ \pi _{K}\circ \iota
_{K}=g\circ id_{K}=g.$

Suppose that $h^{\prime }:K\rightarrow M$ is an $S$-linear map such that $%
f\circ h^{\prime }=g.$ Since $P$ is $M$-$e$-projective and $f\circ
(h^{\prime }\circ \pi _{K})=(f\circ h^{\prime })\circ \pi _{K}=g\circ \pi
_{K}$, there exist $S$-linear maps $h_{1}^{\prime },h_{2}^{\prime
}:P\rightarrow M$ such that $f\circ h_{1}^{\prime }=0=f\circ h_{2}^{\prime }$
and $h^{\ast }+h_{1}^{\prime }=h^{\prime }\circ \pi _{K}+h_{2}^{\prime }$.
Consider $h_{1}:=h_{1}^{\prime }\circ \iota _{K}$ and $h_{2}:=h_{2}^{\prime
}\circ \iota _{K}.$%
\begin{equation*}
\xymatrix{ K \ar[r]^{\iota_K} & P \ar@<.5ex>[r]^{h'_1} \ar@<-.5ex>[r]_{h'_2}
& M \ar[r]^{f} & N \ar[r] & 0\\ & & & K \ar[u]_{g} \ar[lu]_{h}
\ar@/^0.5pc/[lu]^{h'} \ar@<-.5ex>[d]_{\iota_K} \\ & & & P
\ar@<-.5ex>[u]_{\pi_K} \ar@/^1.0pc/[luu]^{h^{*}}}
\end{equation*}%
Then $f\circ h_{1}=f\circ h_{1}^{\prime }\circ \iota _{K}=0,$ $f\circ
h_{2}=f\circ h_{2}^{\prime }\circ \iota _{K}=0$, and%
\begin{equation*}
\begin{array}{ccccc}
h+h_{1} & = & h^{\ast }\circ \iota _{K}+h_{1}^{\prime }\circ \iota _{K} & =
& (h^{\ast }+h_{1}^{\prime })\circ \iota _{K} \\ 
& = & (h^{\prime }\circ \pi _{K}+h_{2}^{\prime })\circ \iota _{K} & = & 
h^{\prime }\circ \pi _{K}\circ \iota _{K}+h_{2}^{\prime }\circ \iota _{K} \\ 
& = & h^{\prime }+h_{2}. &  & 
\end{array}%
\end{equation*}%
Consequently, $K$ is $M$-$e$-projective.$\blacksquare $
\end{Beweis}
\end{lem}

The following result is a relative version of \cite[Corollary 3.3 (3)]%
{AIKN2018}; moreover, we give a \emph{detailed homological} proof as the one 
\cite{AIKN2018} is compact and categorical.

\begin{prop}
\label{dsum-e-proj}Let $\{P_{i}\}_{I\in I}$ be a family of left $S$%
-semimodules and $M$ a left $S$-semimodule. Then $\bigoplus\limits_{i\in
I}P_{i}$ is $M$-$e$-projective if and only if $P_{i}$ is $M$-$e$-projective
for each $i\in I.$ The class of $e$-projective left $S$-semimodules is
closed under direct sums.
\end{prop}

\begin{Beweis}
($\Longrightarrow $)\ This implication follows by Lemma \ref{ret-proj}.

($\Longleftarrow $)\ Let $g:M\rightarrow N$ be a normal epimorphism and $%
f:\bigoplus\limits_{i\in I}P_{i}\rightarrow N$ be an $S$-linear map. For
every $j\in I,$ there exists an $S$-linear map $h_{j}:P_{j}\rightarrow M$
such that $f\circ \iota _{j}=g\circ h_{j}$, where $\iota
_{j}:P_{j}\longrightarrow \bigoplus\limits_{i\in I}P_{i}$ is the canonical
embedding.%
\begin{equation*}
\xymatrix{M \ar[r]^{g} & N \ar[r] & 0\\ & {\bigoplus\limits_{i\in I} P_i}
\ar[u]_{f} \ar@{-->}[lu]_{h}\\ & {P_j} \ar[u]_{\iota_j} \ar@{-->}[luu]^{h_j}}
\end{equation*}

By the \emph{Universal Property of Direct Coproducts}, there exists a unique 
$S$-linear map $h:\bigoplus\limits_{i\in I}P_{i}\rightarrow M$ such that $%
h\circ \iota _{j}=h_{j}$ for every $j\in I,$ i.e.%
\begin{equation*}
h:\bigoplus\limits_{i\in I}P_{i}\rightarrow M,\text{ }\sum\limits_{i\in
I}p_{i}\mapsto \sum\limits_{i\in I}h_{i}(p_{i}).
\end{equation*}%
Notice that $h$ is $S$-linear and well defined since the sum $%
\sum\limits_{i\in I}p_{i}$ is finite (all but finitely many of the
coordinates are zero). Moreover, we have%
\begin{equation*}
\begin{array}{ccccc}
(g\circ h)(\sum\limits_{i\in I}p_{i}) & = & g(\sum\limits_{i\in
I}h_{i}(p_{i})) & = & \sum\limits_{i\in I}(g\circ h_{i})(p_{i}) \\ 
& = & \sum\limits_{i\in I}(f\circ \iota _{i})(p_{i}) & = & 
f(\sum\limits_{i\in I}\iota _{i}(p_{i})) \\ 
& = & f(\sum\limits_{i\in I}p_{i}). &  & 
\end{array}%
\end{equation*}

Suppose that $h^{\prime }:\bigoplus\limits_{i\in I}P_{i}\rightarrow M$ is an 
$S$-linear map with $g\circ h^{\prime }=f.$ Then $f\circ \iota _{j}=g\circ
h^{\prime }\iota _{j}$ for every $j\in I.$ Since $P_{j}$ is $e$-projective
for every $j\in I,$ there exist $S$-linear maps $\tilde{h}_{j},\hat{h}%
_{j}:P_{j}\rightarrow M$ such that $g\circ \tilde{h}_{j}=0=g\circ \hat{h}%
_{j} $ and $h_{j}+\tilde{h}_{j}=h_{j}^{\prime }+\hat{h}_{j}$.

By the \emph{Universal Property of Direct Coproducts}, there exist $S$%
-linear maps%
\begin{equation*}
\tilde{h},\hat{h}:\bigoplus\limits_{i\in I}P_{i}\rightarrow M,\text{ }\tilde{%
h}(\sum\limits_{i\in I}p_{i}):=\sum\limits_{i\in I}\tilde{h}_{i}(p_{i})\text{
and }\hat{h}(\sum\limits_{i\in I}p_{i})=\sum\limits_{i\in I}\hat{h}%
_{i}(p_{i}).
\end{equation*}%
\begin{equation*}
\xymatrix{{P} \ar@<-.5ex>[r]_{\tilde{h}} \ar@<.5ex>[r]^{\hat{h}} & M
\ar[r]^{g} & N \ar[r] & 0\\ & & {\bigoplus\limits_{i\in I} P_i} \ar[u]_{f}
\ar[lu]^{h} \ar@/^-0.5pc/[lu]_{h'}\\ & & {P_j} \ar[u]_{\iota_j}
\ar@/^0.5pc/[luu]^{h_j}}
\end{equation*}%
Both maps are $S$-linear, and well defined since the sum $\sum\limits_{i\in
I}p_{i}$ is finite (all but finitely many of the coordinates are zero).
Moreover, we have 
\begin{eqnarray*}
(g\circ \tilde{h})(\sum\limits_{i\in I}p_{i}) &=&g(\sum\limits_{i\in I}%
\tilde{h}_{i}(p_{i}))=\sum\limits_{i\in I}(g\circ \tilde{h}_{i})(p_{i})=0; \\
\hat{h}(\sum\limits_{i\in I}p_{i}) &=&g(\sum\limits_{i\in I}\hat{h}%
_{i}(p_{i}))=\sum\limits_{i\in I}(g\circ \hat{h}_{i})(p_{i})=0
\end{eqnarray*}%
and%
\begin{equation*}
\begin{array}{ccccc}
(h+\tilde{h})(\sum\limits_{i\in I}p_{i}) & = & h(\sum\limits_{i\in I}p_{i})+%
\tilde{h}(\sum\limits_{i\in I}p_{i}) & = & \sum\limits_{i\in
I}h_{i}(p_{i})+\sum\limits_{i\in I}\tilde{h}_{i}(p_{i}) \\ 
& = & \sum\limits_{i\in I}(h_{i}+\tilde{h}_{i})(p_{i}) & = & 
\sum\limits_{i\in I}(h_{i}^{\prime }+\hat{h}_{i})(p_{i}) \\ 
& = & (h^{\prime }+\hat{h})(\sum\limits_{i\in I}p_{i}). &  & 
\end{array}%
\end{equation*}%
Hence $\bigoplus\limits_{i\in I}P_{i}$ is $M$-$e$-projective. $\blacksquare $
\end{Beweis}

\begin{prop}
\label{lem182}Let $P$ be a left $S$-semimodule. If%
\begin{equation*}
0\longrightarrow K\overset{\iota }{\longrightarrow }L\overset{\pi }{%
\longrightarrow }M\longrightarrow 0
\end{equation*}%
is an exact sequence of left $S$-semimodules and $P$ is $L$-$e$-projective,
then $P$ is $K$-$e$-projective and $M$-$e$-projective.
\end{prop}

\begin{lem}
\begin{Beweis}
Assume that $P$ is $L$-$e$-projective.

\begin{itemize}
\item \textbf{Claim I:\ }$P$ is $M$-$e$-projective. Let $f:M\rightarrow N$
be a normal epimorphism and $g:P\rightarrow N$ an $S$-linear map.%
\begin{equation*}
\xymatrix{ L \ar[r]^{\pi} & M \ar[r]^{f} & N \ar[r] & 0 \\ & & P \ar[u]_{g}
\ar@{-->}[llu]^{h}}
\end{equation*}%
Since $\pi $ and $f$ are normal epimorphism, $f\circ \pi $ is a normal
epimorphism as well (by Lemma \ref{i-normal} (2)(c)). Since $P$ is $L$-$e$%
-projective, there exists an $S$-linear map $h:P\rightarrow M$ such that $%
f\circ \pi \circ h=g.$ Then $\pi \circ h:P\rightarrow M$ is an $S$-linear
map satisfying $f\circ (\pi \circ h)=g.$

Suppose there exists an $S$-linear map $h^{\prime }:P\rightarrow M$ such
that $f\circ h^{\prime }=g.$ Since $\pi $ is a normal epimorphism and $P$ is 
$L$-$e$-projective, there exists an $S$-linear map $h^{\ast }:P\rightarrow L$
such that $\pi \circ h^{\ast }=h^{\prime }.$%
\begin{equation*}
\xymatrix{ L \ar[r]^{\pi} & M \ar[r] & 0\\ & P \ar[u]_{h'}
\ar@{-->}[lu]^{h^*}}
\end{equation*}%
Moreover, $(f\circ \pi )\circ h^{\ast }=f\circ (\pi \circ h^{\ast })=f\circ
h^{\prime }=g.$ Since $P$ is $L$-$e$-projective, there exist $S$-linear maps 
$h_{1},h_{2}:P\rightarrow L$ such that $f\circ \pi \circ h_{1}=0=f\circ \pi
\circ h_{2}$ and $h+h_{1}=h^{\ast }+h_{2}.$ 
\begin{equation*}
\xymatrix{ P \ar@<.5ex>[rr]^{h_2} \ar@<-.5ex>[rr]_{h_1} & & L \ar[rr]^{\pi}
& & M \ar[rr]^{f} & & N \\ \\ & & & & & & P \ar[uu]_{g} \ar[lllluu]_{h^*}
\ar@/^0.5pc/[lllluu]^{h} \ar[lluu]_{h'} }
\end{equation*}%
Thus, $\pi \circ h_{1},$ $\pi \circ h_{2}:P\rightarrow M$ are $S$-linear
maps such that $f\circ \pi \circ h_{1}=0=f\circ \pi \circ h_{2}.$ Moreover, 
\begin{equation*}
\pi \circ h+\pi \circ h_{1}=\pi \circ (h+h_{1})=\pi \circ (h^{\ast
}+h_{2})=\pi \circ h^{\ast }+\pi \circ h_{2}=h^{\prime }+\pi \circ h_{2}.
\end{equation*}%
Consequently, $P$ is $M$-$e$-projective.

\item \textbf{Claim II:\ }$P$ is $K$-$e$-projective. Let $f:K\rightarrow N$
be a normal $S$-epimorphism and $g:P\rightarrow N$ an $S$-linear map. By
Corollary \ref{P-con}, $(\iota ^{\prime },f^{\prime };Q:=(N\oplus L)/\rho )$
is a pushout of $(f,\iota )$ such that $\rho $ is a congruence relation on $%
N\oplus L$ and 
\begin{equation*}
\iota ^{\prime }:N\longrightarrow Q,\text{ }n\mapsto \lbrack (n,0)]_{\rho }%
\text{ and }f^{\prime }:L\longrightarrow Q,\text{ }l\mapsto \lbrack
(0,l)]_{\rho }.
\end{equation*}%
\begin{equation*}
\xymatrix{ L \ar[r]^{f'} & Q \\ K \ar[r]^{f} \ar[u]^{\iota} & N \ar[r]
\ar[u]_{\iota'} & 0 \\ & P \ar[u]_{g} \ar@/^3.0pc/@{-->}[luu]^{h} }
\end{equation*}%
Since $\iota $ is a normal $S$-monomorphism and $f$ is a normal $S$%
-epimorphism, it follows by Lemma \ref{transfers} (2) $\&\ $(4) that $\iota
^{\prime }$ is a normal monomorphism and it follows, by Lemma \ref{transfers}
(3), that $f^{\prime }$ is a normal epimorphism. Since $f^{\prime }$ is a
normal epimorphism and $P$ is $L$-$e$-projective, there exists an $S$-linear
map $h:P\rightarrow L$ such that $f^{\prime }\circ h=\iota ^{\prime }\circ g$%
.

Let $p\in P.$ Since $f$ is surjective, there exists $k\in K$ such that $%
f(k)=g(p).$ Notice that $(f^{\prime }\circ \iota )(k)=(\iota ^{\prime }\circ
f)(k)=(\iota ^{\prime }\circ g)(p)=(f^{\prime }\circ h)(p).$ Since $%
f^{\prime }$ is $k$-normal, there exist $l_{1},l_{2}\in Ker(f^{\prime })$
such that $\iota (k)+l_{1}=h(p)+l_{2}.$

Let 
\begin{equation*}
CP=(\iota ^{\ast },f^{\ast };(N\oplus L)/\rho ^{\ast })
\end{equation*}%
be the $C$-pushout of $(f,\iota )$ (defined in \ref{cpush}). Since $Q$ is a
pushout of $(f,\iota ),$ there exists, by the\emph{\ Universal Property of
Pushouts}, an $S$-linear map $\varphi :Q\rightarrow (N\oplus L)/\rho ^{\ast
} $ such that $\varphi \circ \iota ^{\prime }=\iota ^{\ast }$ and $\varphi
\circ f^{\prime }=f^{\ast }$. Notice that for $i=1,2:$%
\begin{equation*}
\lbrack (0,l_{i})]_{\rho ^{\ast }}=f^{\ast }(l_{i})=\varphi \circ f^{\prime
}(l_{i})=\varphi (0)=[(0,0)]_{\rho ^{\ast }},
\end{equation*}%
and so there exist $k_{i_{1}},k_{i_{2}}\in K$ such that $%
f(k_{i_{1}})=f(k_{i_{2}})$ and $l_{i}+\iota (k_{i_{2}})=\iota (k_{i_{1}})$.

Since $\iota $ is a normal monomorphism, $\iota (K)\subseteq L$ is
subtractive, whence $l_{1},l_{2}\in \iota (K)$, i.e. $l_{1}=\iota (k_{1})$
and $l_{2}=\iota (k_{2})$ for some $k_{1},k_{2}\in K.$ It follows that $%
\iota (k)+\iota (k_{1})=h(p)+\iota (k_{2}),$ we conclude that $h(p)\in \iota
(K)$ (as $\iota $ is a normal monomorphism). Let $k_{p}\in K$ be such that $%
h(p)=\iota (k_{p}).$ Notice that this $k_{p}$ is unique since $\iota $ is an
injective. Therefore%
\begin{equation*}
h^{\prime }:P\longrightarrow K,\text{ }p\mapsto k_{p}
\end{equation*}%
is well defined. Clearly, $h^{\prime }$ is $S$-linear. Now, for every $p\in
P,$ we have 
\begin{equation*}
(\iota ^{\prime }\circ f\circ h^{\prime })(p)=(f^{\prime }\circ (\iota \circ
h^{\prime }))(p)=(f^{\prime }\circ \iota )(k_{p})=(f^{\prime }\circ
h)(p)=(\iota ^{\prime }\circ g)(p),
\end{equation*}
whence $(f\circ h^{\prime })(p)=g(p)$ as $\iota ^{\prime }$ is injective.

Suppose that there exists an $S$-linear map $h^{\ast }:P\rightarrow K$ such
that $f\circ h^{\ast }=g.$ It follows that $f^{\prime }\circ \iota \circ
h^{\ast }=\iota ^{\prime }\circ f\circ h^{\ast }=\iota ^{\prime }\circ g.$
Since $P$ is $L$-$e$-projective, there exist $S$-linear maps $%
h_{1},h_{2}:P\rightarrow L$ such that $f^{\prime }\circ h_{1}=0=f^{\prime
}\circ h_{2}$ and $h+h_{1}=\iota \circ h^{\ast }+h_{2}.$ Let $p\in P.$ For $%
i=1,2,$ and since $h_{i}(p)\in Ker(f^{\prime }),$ there exists $k_{p}^{i}\in
K$ such that $h_{i}(p)=\iota (k_{p}^{i})$ (which is indeed unique as $\iota $
injective). Then we have two well defined maps%
\begin{equation*}
h_{1}^{\prime }:P\longrightarrow K,\text{ }p\mapsto k_{p}^{1}\text{ and }%
h_{2}^{\prime }:P\longrightarrow K,\text{ }p\mapsto k_{p}^{2}.
\end{equation*}%
which can be easily shown to be $S$-linear.

For every $p\in P,$ and for $i=1,2$ we have 
\begin{equation*}
(\iota ^{\prime }\circ f\circ h_{i}^{\prime })(p)=(\iota ^{\prime }\circ
f)(k_{p}^{i})=(f^{\prime }\circ \iota )(k_{p}^{i})=(f^{\prime }\circ
h_{i})(p)=0,
\end{equation*}%
whence $(f\circ h_{i}^{\prime })(p)=0$ as $\iota ^{\prime }$ is injective.
Moreover, we have%
\begin{equation*}
\begin{array}{ccccc}
\iota ((h^{\ast }+h_{2}^{\prime })(p)) & = & (\iota \circ h^{\ast
})(p)+(\iota \circ h_{2}^{\prime })(p) & = & (\iota \circ h^{\ast
})(p)+\iota (k_{p}^{2}) \\ 
& = & (\iota \circ h^{\ast })(p)+h_{2}(p) & = & (\iota \circ h^{\ast
}+h_{2})(p) \\ 
& = & (h+h_{1})(p) & = & h(p)+h_{1}(p) \\ 
& = & \iota (k_{p})+\iota (k_{p}^{1}) & = & (\iota \circ h^{\prime
})(p)+(\iota \circ h_{1}^{\prime })(p) \\ 
& = & \iota ((h^{\prime }+h_{1}^{\prime })(p)) &  & 
\end{array}%
\end{equation*}%
whence $(h^{\ast }+h_{2}^{\prime })(p)=(h^{\prime }+h_{1}^{\prime })(p)$ as $%
\iota $ is injective. Consequently, $P$ is $K$-$e$-projective.$\blacksquare $
\end{itemize}
\end{Beweis}
\end{lem}

\begin{ex}
\label{nonseproj}Consider the semiring $S:=M_{2}(\mathbb{R}^{+})$ and the
subtractive left ideal%
\begin{equation*}
N_{1}=\left\{ \left[ {%
\begin{array}{cc}
a & a \\ 
b & b%
\end{array}%
}\right] |\text{ }a,b\in \mathbb{R}^{+}\right\} .
\end{equation*}%
Then $S/N_{1}$ is not an $S$-$k$-projective $S$-semimodule, whence not $S$-$%
e $-projective.
\end{ex}

\begin{Beweis}
Let $\pi :S\rightarrow S/N_{1}$ be the canonical map and $id_{S/N_{1}}$ be
the identity map of $S/N_{1}$. Notice that $\pi $ is a normal epimorphism.
Consider 
\begin{equation*}
e_{1}=\left[ {%
\begin{array}{cc}
1 & 0 \\ 
0 & 0%
\end{array}%
}\right] \text{ and }e_{2}=\left[ {%
\begin{array}{cc}
0 & 0 \\ 
0 & 1%
\end{array}%
}\right] .
\end{equation*}%
Suppose that there exists an $S$-linear map $g:S/N_{1}\rightarrow S$ such
that $\pi g=id_{S/N_{1}}$. Then $g(\overline{e_{1}})\in \pi ^{-1}(\overline{%
e_{1}})$ and $g(\overline{e_{2}})\in \pi ^{-1}(\overline{e_{2}})$. Write $g(%
\overline{e_{1}})=\left[ {%
\begin{array}{cc}
p & q \\ 
r & s%
\end{array}%
}\right] $ for some $p,q,r,s\in \mathbb{R}^{+}$. Then 
\begin{equation*}
\left[ {%
\begin{array}{cc}
p+k & q+k \\ 
r+l & s+l%
\end{array}%
}\right] =\left[ {%
\begin{array}{cc}
m+1 & m \\ 
n & n%
\end{array}%
}\right]
\end{equation*}%
for some $k,l,m,n\in \mathbb{R}^{+}$, which implies that $r=s$ and $p=q+1$
as $\mathbb{R}^{+}$ is cancellative. By relabeling, we have 
\begin{equation*}
g(\overline{e_{1}})=\left[ {%
\begin{array}{cc}
a+1 & a \\ 
b & b%
\end{array}%
}\right] \text{ and }g(\overline{e_{2}})=\left[ {%
\begin{array}{cc}
c & c \\ 
d & d+1%
\end{array}%
}\right] \text{ for some }a,b,c,d\in \mathbb{R}^{+}.
\end{equation*}

Let $x:=\left[ {%
\begin{array}{cc}
p & q \\ 
r & s%
\end{array}%
}\right] \in S.$ Then 
\begin{equation*}
x=\left[ {%
\begin{array}{cc}
p & 0 \\ 
r & 0%
\end{array}%
}\right] e_{1}+\left[ {%
\begin{array}{cc}
0 & q \\ 
0 & s%
\end{array}%
}\right] e_{2},
\end{equation*}%
which implies that 
\begin{equation*}
g(\overline{x})=\left[ {%
\begin{array}{cc}
p & 0 \\ 
r & 0%
\end{array}%
}\right] g(\overline{e_{1}})+\left[ {%
\begin{array}{cc}
0 & q \\ 
0 & s%
\end{array}%
}\right] g(\overline{e_{2}})=\left[ {%
\begin{array}{cc}
pa+dq+p & pa+dq+q \\ 
ra+sd+r & ra+sd+s%
\end{array}%
}\right] .
\end{equation*}%
But%
\begin{equation*}
x=\left[ {%
\begin{array}{cc}
p & 1 \\ 
r & 0%
\end{array}%
}\right] e_{1}+\left[ {%
\begin{array}{cc}
0 & q \\ 
1 & s%
\end{array}%
}\right] e_{2},
\end{equation*}%
which implies%
\begin{equation*}
\begin{array}{ccc}
\left[ {%
\begin{array}{cc}
pa+dq+p & pa+dq+q \\ 
ra+sd+r & ra+sd+s%
\end{array}%
}\right] & = & g(\overline{x})=\left[ {%
\begin{array}{cc}
p & 1 \\ 
r & 0%
\end{array}%
}\right] g(\overline{e_{1}})+\left[ {%
\begin{array}{cc}
0 & q \\ 
1 & s%
\end{array}%
}\right] g(\overline{e_{2}}) \\ 
& = & \left[ {%
\begin{array}{cc}
(pa+dq+p)+b & (pa+dq+q)+b \\ 
(ra+sd+r)+c & (ra+sd+s)+c%
\end{array}%
}\right] ,%
\end{array}%
\end{equation*}%
whence $b=0=c$ as $\mathbb{R}^{+}$ is cancellative. Thus 
\begin{equation*}
g(\overline{e_{1}})=\left[ {%
\begin{array}{cc}
a+1 & a \\ 
0 & 0%
\end{array}%
}\right] \text{ and }g(\overline{e_{2}})=\left[ {%
\begin{array}{cc}
0 & 0 \\ 
d & d+1%
\end{array}%
}\right] \text{ for some }a,d\in \mathbb{R}^{+}.
\end{equation*}

Let $y:=\left[ {%
\begin{array}{cc}
2 & 1 \\ 
0 & 0%
\end{array}%
}\right] $. Notice that $\overline{e_{1}}=\overline{y}$, whence 
\begin{equation*}
\left[ {%
\begin{array}{cc}
a+1 & a \\ 
0 & 0%
\end{array}%
}\right] =g(\overline{e_{1}})=g(\overline{y})=\left[ {%
\begin{array}{cc}
2a+d+2 & 2a+d+1 \\ 
0 & 0%
\end{array}%
}\right] ,
\end{equation*}%
and so $a=2a+d+1$. Since $\mathbb{R}^{+}$ is cancellative, $a+d+1=0$, that
is $1$ has additive inverse, a contradiction. Hence, there is no such $S$%
-linear map $g$ with $\pi \circ g=id_{S/I},$ \emph{i.e.}, $S/I$ is not $S$-$%
k $-projective. Since $S/I$ is not $S$-$k$-projective, $S/I$ is not $S$-$e$%
-projective.$\blacksquare $
\end{Beweis}

Recall the following fact about the relative projectivity for modules over
rings.

\begin{punto}
\label{Opl-proj}\cite{Wis1991} Let $R$ be a ring, $P$ a left $R$-module and $%
\{M_{\lambda }\}_{\lambda \in \Lambda }$ a collection of left $S$%
-semimodules such that $P$ is $M_{\lambda }$-projective for every $\lambda
\in \Lambda .$ If $\Lambda =\{\lambda _{1},\cdots ,\lambda _{k}\}$ is
finite, then $P$ is $\bigoplus\limits_{n=1}^{k}M_{\lambda _{n}}$-projective.
If $_{R}P$ is finitely generated and $\Lambda $ is arbitrary, then $P$ is $%
\bigoplus\limits_{\lambda \in \Lambda }M_{\lambda }$-projective (even if $%
\Lambda $ is infinite).
\end{punto}

We provide a counter example showing that the result corresponding to \ref%
{Opl-proj} for the relative $e$-projectivity for semimodules over semiring
does \emph{not} necessarily hold. The same example serves to show that the
converse of Proposition \ref{lem182} is not true (even when $M=L\oplus N$).

\begin{ex}
\label{endsprojnotmid}Let $S:=M_{2}(\mathbb{R}^{+})=E_{1}\oplus E_{2},$ where%
\begin{equation*}
E_{1}=\left\{ \left[ {%
\begin{array}{cc}
a & 0 \\ 
b & 0%
\end{array}%
}\right] |\text{ }a,b\in \mathbb{R}^{+}\right\} ,\text{ }E_{2}:=\left\{ %
\left[ {%
\begin{array}{cc}
0 & c \\ 
0 & d%
\end{array}%
}\right] |\text{ }c,d\in \mathbb{R}^{+}\right\}
\end{equation*}%
and consider%
\begin{equation*}
K=\left\{ \left[ {%
\begin{array}{cc}
u & u \\ 
v & v%
\end{array}%
}\right] |\text{ }u,v\in \mathbb{R}^{+}\right\} \text{ and }P:=S/K.
\end{equation*}%
Then 
\begin{equation*}
0\rightarrow E_{1}\xrightarrow{\iota_{E_1}}S\xrightarrow{\pi_{E_2}}%
E_{2}\rightarrow 0
\end{equation*}%
is exact, $P$ is $E_{1}$-$e$-projective and $E_{2}$-$e$-projective. However, 
$P$ is not $S$-$e$-projective (notice that $S=E_{1}\oplus E_{2}$).
\end{ex}

\begin{Beweis}
Since $E_{1}\oplus E_{2}=S$, it follows by the proof of Example \ref%
{nonseproj} that $P$ is not $(E_{1}\oplus E_{2})$-$e$-projective. Notice
that $E_{1}$ and $E_{2}$ are ideal-simple left $S$-subsemimodules of $S.$
Let $L\neq 0$ and $f:E_{1}\rightarrow L$ be a normal $S$-epimorphism. Then $%
Ker(f)\subsetneqq E_{1},$ whence $Ker(f)=0$ as $E_{1}$ is ideal-simple.
Since $f$ is $k$-normal and $Ker(f)=0,$ $f$ is injective, whence an
isomorphism. If $g:P\rightarrow L$ is an $S$-linear map, then $f^{-1}\circ
g:P\rightarrow E_{1}$ is an $S$-linear map such that $f\circ f^{-1}\circ
g=g, $ and whenever there exists an $S$-linear map $h:P\rightarrow E_{1}$
such that $f\circ h=g,$ we have $h=f^{-1}\circ f\circ h=f^{-1}\circ g.$
Hence, $P$ is $E_{1}$-$e$-projective. Similarly, one can prove that $P$ is $%
E_{2}$-$e$-projective.$\blacksquare $
\end{Beweis}


\end{document}